\documentclass{article}

\usepackage{caption}
\usepackage{subcaption}
\usepackage{arxiv}
\usepackage{amsmath}
\usepackage{amssymb}
\usepackage{mathtools}
\usepackage{amsthm}
\usepackage[utf8]{inputenc} 
\usepackage[T1]{fontenc}    
\usepackage{hyperref}       
\usepackage{url}            
\usepackage{booktabs}       
\usepackage{amsfonts}       
\usepackage{nicefrac}       
\usepackage{microtype} 

\usepackage{hyperref}
\PassOptionsToPackage{compress, round}{natbib}
\usepackage{natbib}

\usepackage[capitalize,noabbrev]{cleveref}

\usepackage{lipsum}
\usepackage{graphicx}

\newcommand{\E}{\mathbb{E}}

\newcommand{\argmax}{\text{argmax}}
\newcommand{\argmin}{\text{argmin}}
\newcommand{\ceil}[1]{\left\lceil#1\right\rceil}

\newcommand{\tx}{\tilde{\vx}}
\newcommand{\ty}{\tilde{\vy}}
\newcommand{\tz}{\tilde{\vz}}

\newcommand{\re}{\operatorname{Re}}
\newcommand{\im}{\operatorname{Im}}

\newcommand{\bigotilde}{\tilde{\mathcal{O}}}

\newcommand{\cO}{\mathcal{O}}

\newcommand{\cE}{\mathcal{E}}

\newcommand{\cI}{\mathcal{I}}

\newcommand{\cU}{\mathcal{U}}
\newcommand{\cX}{\mathcal{X}}
\newcommand{\cY}{\mathcal{Y}}

\newcommand{\norm}[1]{\left\|#1\right\|}
\newcommand{\abs}[1]{\left\vert#1\right\vert}

\newcommand{\R}{\mathbb{R}}

\newcommand{\vx}{\boldsymbol{x}}
\newcommand{\vy}{\boldsymbol{y}}

\newcommand{\vz}{\boldsymbol{z}}

\DeclareMathAlphabet{\mathsfit}{T1}{\sfdefault}{\mddefault}{\sldefault}
\SetMathAlphabet{\mathsfit}{bold}{T1}{\sfdefault}{\bfdefault}{\sldefault}

\makeatletter
\newcommand{\vast}{\bBigg@{3}}
\newcommand{\Vast}{\bBigg@{3.5}}
\makeatother

\newcommand{\haochuan}[1]{\textcolor{blue}{\{{\textbf{HL}: \em #1}\}}}

\newcommand{\tr}{\text{Trace}}

\newcommand{\C}{\mathbb{C}}

\newcommand{\ky}{\kappa}
\newcommand{\kx}{\kappa_{\vx}}

\newcommand{\etax}{\eta_{\vx}}
\newcommand{\etay}{\eta_{\vy}}
\newcommand{\mux}{\mu_{\vx}}

\newcommand{\vxi}{\boldsymbol{\xi}}
\newcommand{\0}{\mathbf{0}}

\theoremstyle{plain}
\newtheorem{theorem}{Theorem}[section]

\newtheorem{lemma}[theorem]{Lemma}
\newtheorem{corollary}[theorem]{Corollary}
\theoremstyle{definition}
\newtheorem{definition}[theorem]{Definition}
\newtheorem{assumption}[theorem]{Assumption}
\theoremstyle{remark}
\newtheorem{remark}[theorem]{Remark}

\title{On Convergence of Gradient Descent Ascent: A Tight Local Analysis}

\author{
 Haochuan Li \\
  MIT EECS\\
  \texttt{haochuan@mit.edu} \\
   \And
 Farzan Farnia \\
  CUHK CSE\\
  \texttt{farnia@cse.cuhk.edu.hk} \\
  \And
 Subhro Das \\
  MIT-IBM AI Lab, IBM Research\\
  \texttt{Subhro.Das@ibm.com } \\
    \And
 Ali Jadbabaie \\
  MIT CEE\\
  \texttt{jadbabai@mit.edu} 
}

\begin{document}
\maketitle

\begin{abstract}
Gradient Descent Ascent (GDA) methods are the mainstream algorithms for minimax optimization in generative adversarial networks (GANs). Convergence properties of GDA have drawn significant interest in the recent literature. Specifically, for $\min_{\vx} \max_{\vy} f(\vx;\vy)$ where $f$ is strongly-concave in $\vy$ and possibly nonconvex in $\vx$, (Lin et al., 2020) proved the convergence of GDA with a stepsize ratio $\etay/\etax=\Theta(\kappa^2)$ where $\etax$ and $\etay$ are the stepsizes for $\vx$ and $\vy$ and $\kappa$ is the condition number for $\vy$. While this stepsize ratio suggests a slow training of the min player, practical GAN algorithms typically adopt similar stepsizes for both variables, indicating a wide gap between theoretical and empirical results. In this paper, we aim to bridge this gap by analyzing the \emph{local convergence} of general \emph{nonconvex-nonconcave} minimax problems. We demonstrate that a stepsize ratio of $\Theta(\kappa)$ is necessary and sufficient for local convergence of GDA to a Stackelberg Equilibrium, where $\kappa$ is the local condition number for $\vy$. We prove a nearly tight convergence rate with a matching lower bound. We further extend the convergence guarantees to stochastic GDA and extra-gradient methods (EG). Finally, we conduct several numerical experiments to support our theoretical findings.
\end{abstract}

\section{Introduction}
\label{sec:intro}

Minimax learning frameworks including generative adversarial networks (GANs) \citep{goodfellow2014generative} and adversarial training \citep{madry2017towards} have achieved great success in various machine learning tasks. According to these frameworks, the underlying machine learning problem is formulated as a zero-sum game between two min and max players trying to optimize a learning objective in the opposite directions. Therefore, minimax learning tasks are commonly formulated through a minimax optimization problem of the following form:
\begin{align}\label{eq:minimax}
\min_{\vx\in\cX}\max_{\vy\in\cY}f(\vx;\vy).
\end{align}
In the above problem, $\vx$ and $\vy$ denote the optimization variables for the minimization and maximization subproblems, respectively, and $f(\vx;\vy)$ represents the minimax learning objective estimated over observed training data. For example, in GAN settings, $\vx$ contains the parameters of a generator machine producing real-like samples, and $\vy$ represents the parameters of a discriminator machine distinguishing the generated samples from real training data. 

Since the machine players of modern minimax learning frameworks are typically chosen as deep neural networks, these frameworks lead to difficult nonconvex-nonconcave minimax optimization problems where $f(\vx;\vy)$ is possibly nonconvex in $\vx$ and nonconcave in $\vy$. In such optimization problems, standard gradient-based algorithms are not guaranteed to find a locally optimal solution. However, numerous empirical studies on GANs suggest that a simultaneous-update gradient-based optimization algorithm such as Gradient Descent Ascent (GDA) can successfully train a satisfactory generator function. The wide gap between standard convergence guarantees and empirical results in GAN experiments has inspired several recent studies on the convergence behavior of GDA-type methods in general nonconvex-nonconcave minimax problems.

A recent line of works in the minimax optimization literature focuses on the convergence properties of gradient-based optimization algorithms in the \textit{nonconvex-concave minimax settings} where the objective $f(\vx;\vy)$ is further assumed to be a concave function of $\vy$. Regarding this class of minimax optimization problems, a well-known result proved by \cite{lin2020gradient} demonstrates a convergence guarantee for the GDA and stochastic GDA (SGDA) algorithms in the smooth nonconvex-strongly-concave (NC-SC) minimax settings where $f(\vx;\vy)$ is assumed to be $L$-smooth and $\mu$-strongly-concave in $\vy$. Under these assumptions, \citet{lin2020gradient} analyze the convergence of two-time-scale GDA with different stepsizes $\etax=\Theta(\frac{1}{\ky^2L})$ and $\etay=\Theta(\frac{1}{L})$ which achieves a gradient complexity bound of $\cO\left(\kappa^2/\epsilon^2\right)$ for $\kappa:=L/\mu$ defined as the condition number of the minimax objective and $\epsilon$ denoting the stationarity degree of the desired solution. 

\begin{figure*}[t]
     \centering
     \begin{subfigure}[b]{0.47\textwidth}
         \centering
         \includegraphics[height=0.47\textwidth]{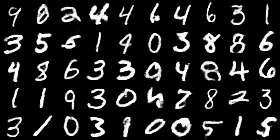}
         \caption{}
         \label{fig:intro_a}
     \end{subfigure}
     \hfill
     \begin{subfigure}[b]{0.47\textwidth}
         \centering
         \includegraphics[height=0.47\textwidth]{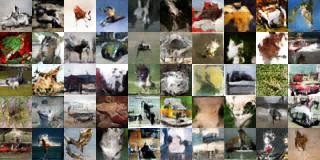}
         \caption{}
         \label{fig:intro_b}
     \end{subfigure}
        \caption{Generated images of the learned generator on MNIST (a) and CIFAR10 (b). For both MNIST and CIFAR10, we train WGAN-GP models~\citep{Gulrajani2017ImprovedTO} using simultaneous GDA with $\etax=\etay=0.001$.}
        \label{fig:intro}
\end{figure*}

While the above result sheds light on the convergence properties of two-time-scale GDA methods for general minimax settings, it suggests choosing a significantly smaller stepsize value for the minimization subproblem compared to that of the maximization subproblem. While the recommended stepsize ratio $\eta_{\vy}/\eta_{\vx}=\Theta(\kappa^2)$ would tend to infinity for general nonconvex-nonconcave problems, practical GAN experiments often choose similar stepsize values for the minimization and maximization tasks. 
For example, Figure~\ref{fig:intro} shows that choosing a $\Theta(1)$ stepsize ratio may already suffice to converge to a good solution for GANs trained over standard MNIST and CIFAR-10 datasets. In addition, \citet{farnia2021train} provide theoretical and numerical evidence that a smaller GDA stepsize ratio can result in better generalization performance for the gradient-based minimax learners.   
Such empirical and theoretical results question the necessity of a large stepsize ratio for minimax learning algorithms, and further motivate the following question: 

\textit{What is the best stepsize ratio that ensures convergence of GDA and what is the corresponding convergence rate? }

To answer this question, we develop a tight local convergence analysis near a Stackelberg Equilibrium of smooth minimax optimization problems. Specifically, we show that a stepsize ratio of $\Theta(\kappa)$, where $\kappa$ is the local condition number of the max player, is necessary and sufficient to ensure the convergence of GDA. When applied to the NC-SC setting, it is much smaller than that in \citep{lin2020gradient} and enables a faster local convergence rate. We also show that our stepsize ratio and the corresponding convergence rate are nearly optimal by providing matching lower bounds.

Different from standard analysis methods in optimization, we leverage the tools in linear dynamical systems and control theory to derive necessary and sufficient conditions under which GDA is guaranteed to converge. Specifically, note that a smooth function is approximately quadratic NC-SC locally around its Stackelberg Equilibrium and thus the dynamics of GDA is nearly a linear time invariant (LTI) dynamical system. The convergence is then determined by the spectral radius $\rho$, the largest modulus of the complex eigenvalues, of the transition matrix of this system. For example, $\rho<1$ implies convergence and $\rho>1$ implies divergence. However, the transition matrix is a non-symmetric and pretty general block matrix whose eigenvalues could be complex numbers and have no closed-form solution. In Lemma~\ref{lem:spectral_M}, we provide a novel and fine-grained analysis of its spectral properties, which is then used to obtain our convergence results. Lemma~\ref{lem:spectral_M} is also our main technical contribution which in itself can be interesting to the optimization community.


\subsection{Summary of main contributions}
\begin{itemize}
    \item  We provide a precise answer to the above question for local convergence to a Stackelberg Equilibrium of locally $L$ smooth functions. In this setting, we prove that choosing $\eta_{\vx}=\Theta(1/r L)$ and $\eta_{\vy}=\Theta(1/L)$ with a stepsize ratio $r=\Omega(\kappa)$ suffices to guarantee convergence.  We then show the necessity of the $\Omega(\kappa)$ stepsize ratio to guarantee convergence. Therefore the optimal choice is $r=\Theta(\kappa)$.
    
    \item For any $r=\Omega(\kappa)$, we provide a nearly tight convergence rate with matching lower bounds in this setting. Our result is locally much faster than that under the conservative stepsize choice in \citep{lin2020gradient}. 
    
    \item
    We further extend our local convergence results to Stochastic Gradient Descent Ascent (SGDA) and extre-gradient methods (EG). For EG, the convergence rate is also optimal for any $r=\Omega(\kappa)$.

    \item We conduct experiments on quadratic NC-SC functions to illustrate the convergence behaviors under different stepsize ratios. We also run the GDA algorithm to optimize non-quadratic functions and show that such behaviors further exist globally for several non-quadratic functions. Finally, we conduct experiments on GANs to show that simultaneous GDA with $\etax\approx\etay$ enjoys fast convergence to a desired solution. 
\end{itemize}

\section{Related work}
\label{sec:related}

\textbf{GDA-based minimax optimization algorithms.} Gradient Descent Ascent (GDA) naturally extends gradient descent (GD) to the minimax optimization problem~\eqref{eq:minimax}. The algorithm takes one step of gradient descent over the minimization variable $\vx$ with stepsize $\eta_{\vx}$ and another step of gradient ascent over the maximization variable $\vy$ with stepsize $\eta_{\vy}$. A lot of convergence results have been established for the average iterates generated by GDA when the objective function is convex in $\vx$ and concave in $\vy$~\citep{Chen1997ConvergenceRI,Nedi2009SubgradientMF,nemirovski2004prox,Du2019LinearCO}.

However, in terms of last iterate convergence, some works \citep{Benam1999MixedEA,Hommes2012MultipleEA,Mertikopoulos2018CyclesIA} show that GDA with equal stepsizes ($\eta_{\vx}=\eta_{\vy}$) can converge to limit cycles or diverge even for convex-concave functions. Two-time-scale GDA with unequal stepsizes ($\eta_{\vx}\neq\eta_{\vy}$) is shown empirically	to alleviate the issues of limit cycles and its local asymptotic convergence to Nash equilibria is theoretically guaranteed~\citep{Heusel2017GANsTB}. \citet{lin2020gradient} established non-asymptotic convergence guarantees in the nonconvex-concave setting.

In a closely-related work, \citet{fiez2021local} demonstrate how to construct a ratio $r^\ast$ that ensures GDA's local convergence. However, this work does not focus on answering what the optimal choice of $r^\ast$ will be. Moreover, the questions of whether the proposed $r^\ast$ is the optimal choice and how large $r^\ast$ and their convergence rates are remain unanswered in this paper. 
In our work, we provide a precise answer by showing what the optimal $r^\ast$ and the convergence rates are for local convergence.

\textbf{Other algorithms.}  Many works study algorithms other than GDA and provide strong theoretical results in the convex-concave setting. The popular extra-gradient method (EG), which takes two gradient updates per iteration, is one of them~\citep{korpelevich1976extragradient, tseng1995linear, nemirovski2004prox, chambolle2011first, yadav2017stabilizing}. Another line studies the optimistic gradient method~\citep{rakhlin2013optimization, daskalakis2017training, gidel2018variational, mertikopoulos2018optimistic, hsieh2019convergence}. \citet{mokhtari2020convergence} provided a unified analysis of both methods as approximations of the classical proximal point method. \cite{lin2020nearoptimal} focused on accelerating the known rates in terms of the conditional number dependence and their rates match the lower bounds~\citep{ouyang2019lower,Ibrahim2019LowerBA,Zhang2019OnLI}.

Recently, several convergence results under the nonconvex-concave setup are established. \citet{rafique2018non} proposed proximally guided gradient methods and analyzed its convergence. \citet{thekumparampil2019efficient} proposed an inexact proximal point method for the nonconvex-concave setting. Many other works also studied different variations of convergence~\citep{nouiehed2019solving, kong2019accelerated, lu2020hybrid, ostrovskii2020efficient}. The algorithm of \citep{lin2020nearoptimal} was shown to be nearly optimal in the nonconvex-strongly-concave setting later in \citep{Li2021ComplexityLB,zhang2021complexity,Han2021LowerCB}

The general nonconvex-nonconcave setting is even more challenging. First, a saddle-point may not exist \citep{jin2020local,farnia2020gans}. Moreover, determining its existence is known to be NP-hard~\citep{daskalakis2020complexity}, and finding an approximate local saddle point is PPAD-complete~\citep{daskalakis2020complexity}. It remains an open problem how to define a nontrivial and tractable notion of minimax points.  

\section{Preliminaries}
\label{sec:pre}
\paragraph{Notation.}
We use bold lower-case letters to denote vectors, upper-case letters to denote matrices, and calligraphic upper-case letters to denote sets. 
Let $\norm{\vx}_2$ be the $\ell_2$ norm of vector $\vx$ and $\norm{A}_2$ be the spectral norm of matrix $A$. 
For a matrix $A$, we use $A^\top$, $\det(A)$, and $\tr(A)$ to denote its transpose, determinant, and trace respectively. Given a symmetric real matrix $A$, we use $\lambda_{\min}(A)$ and $\lambda_{\max}(A)$ to denote its smallest and largest eigenvalues. Given another symmetric real matrix $B$, we write $A< B$ or $B>A$ if $B-A$ is positive definite. We also write $A\le B$ or $B\ge A$ if $B-A$ is positive semi-definite. We use $I_d$ to denote the identity matrix in $\R^{d\times d}$ and omit the subscript when $d$ is clear from the context.
Finally, we use the standard $\cO(\cdot)$, $\Theta(\cdot)$, and $\Omega(\cdot)$ notation, with $\bigotilde(\cdot)$, $\Tilde{\Theta}(\cdot)$, and $\Tilde{\Omega}(\cdot)$ further hiding logarithmic factors.

\subsection{Problem setup}
\label{subsec:setup}
We study the minimax optimization problem \eqref{eq:minimax} where $f$ is twice-continuously differentiable. This assumption guarantees that $\nabla^2 f(\vx;\vy)$ is a symmetric Hessian matrix for any $\vx\in\cX\subseteq \R^n$ and $\vy\in\cY\subseteq\R^m$. We will often denote $\vz=(\vx,\vy)$ for simplicity.

We further assume that there exists a differential Stackelberg Equilibrium $\vz^\ast=(\vx^\ast,\vy^\ast)$. Formally, denote the Hessian at $\vz^\ast$ as
\begin{align}
    \begin{pmatrix}
    \nabla^2_{\vx\vx}f(\vz^\ast) &\nabla^2_{\vx\vy}f(\vz^\ast)\\
    \nabla^2_{\vy\vx}f(\vz^\ast)&\nabla^2_{\vy\vy}f(\vz^\ast)
    \end{pmatrix}=\begin{pmatrix}
    C &B\\
    B^\top&-A
    \end{pmatrix}.
    \label{eq:hessian_def}
\end{align}
Then differential Stackelberg Equilibrium is defined as follows.

\begin{definition}[\cite{pmlr-v119-fiez20a}] A point $\vz^\ast$ satisfying \eqref{eq:hessian_def} is a differential Stackelberg Equilibrium if $\nabla_{\vx} f(\vz^\ast)=\0$, $\nabla_{\vy} f(\vz^\ast)=\0$, $A> 0$,  $C+BA^{-1}B^\top>0$. 
\label{def:stackelberg}
\end{definition}
Here $C+BA^{-1}B^\top$ is the Schur complement with respect to the top left $n\times n$ block in the block matrix $\nabla^2 f(\vz^\ast)$. We will see later that it is also the Hessian of the primal function $\Phi(\cdot)=\max_{\vy\in\cY}f(\cdot;\vy)$ at $\vx^\ast$. The condition $A>0$ implies the function is locally strictly concave in $\vy$ around $\vz^\ast$. The condition $C+BA^{-1}B^\top$ further suggests $\Phi$ is locally strictly convex.

We assume that $f$ is locally $L$ smooth around $z^\ast$. Specifically, assume\footnote{Our definition of smoothness is slightly different from the usual one up to a factor of $2$.}
\begin{align*}
 \max\{\norm{A}_2,\norm{B}_2,\norm{C}_2\}\le L.
\end{align*}
Denote
\begin{align*}
    \mu = \lambda_{\min}(A),\quad\mu_x=\min\{L,\lambda_{\min}(C+BA^{-1}B^\top)\}
\end{align*}
and the corresponding condition numbers 
\begin{align*}
    \ky = L/\mu,\quad \kx=L/\mux.
\end{align*}
We will call $\vz^\ast$ a Stackelberg Equlibrium with parameters $(L,\mu,\mux)$. Also remark that $\Phi$ is locally $(\ky+1)L$ smooth and $\mux$ strongly convex, which implies the local condition number of $\Phi$ is actually $\kx(\ky+1)$.

\paragraph{Learning Dynamics} 
In this paper, we focus on the two-time-scale GDA algorithm and also extend the results to two-time-scale EG. Let $\etax$ and $\etay$ be the stepsizes for $\vx$ and $\vy$ and $r=\etay/\etax$ be the stepsize ratio. 

At each iterate $\vz^k=(\vx^k,\vy^k)$, we are given an exact or stochastic oracle that returns gradient $(g_{\vx}(\vz^k),g_{\vy}(\vz^k))$. If it is an exact oracle, we have
\begin{align*}
    g_{\vx}(\vz^k)=\nabla_{\vx}f(\vz^k),\quad g_{\vy}(\vz^k)=\nabla_{\vy}f(\vz^k).
\end{align*}
In the stochastic setting, we hide the randomness in the random functions $g_{\vx}$ and $g_{\vy}$. Denote the noise of gradient as
\begin{align}
    \vxi_{\vx}^{k}=&g_{\vx}(\vz^k)-\nabla_{\vx}f(\vz^k),\label{eq:gradient_noise}\\\vxi_{\vy}^{k}=&g_{\vy}(\vz^k)-\nabla_{\vy}f(\vz^k).\nonumber
\end{align}
We assume the gradient is estimated with a mini-batch of samples of size $S$ and is unbiased with bounded variance. Formally, we assume
\begin{assumption}
  $\E[\vxi_{\vx}^{k}]=\0$, $\E[\vxi_{\vy}^{k}]=\0$, $\E[\|\vxi_{\vx}^{k}\|_2^2]\le \sigma^2/S$, and $\E[\|\vxi_{\vy}^{k}\|_2^2]\le \sigma^2/S$ for some $\sigma\ge 0$. Also, the noises for different $k$ are independent.
    \label{ass:zero_mean_bounded_variance}
\end{assumption}
Remark that in the stochastic setting, it takes $S$ gradient evaluations to estimate $(g_{\vx}(\vz^k),g_{\vy}(\vz^k))$.

The learning dynamics of (stochastic) GDA is given by
\begin{align}
    \vx^{k+1}=&\;\vx^k-\etax g_{\vx}(\vx^{k};\vy^{k}),\label{eq:gda_def}\\
    \vy^{k+1}=&\;\vy^k+\etay g_{\vy}(\vx^{k};\vy^{k}).\nonumber
\end{align}

The learning dynamics of (stochastic) EG is given by
\begin{align}
    \vx^{k+1/2}=&\;\vx^k-\etax g_{\vx}(\vx^{k};\vy^{k}),\nonumber\\
    \vy^{k+1/2}=&\;\vy^k+\etay g_{\vy}(\vx^{k};\vy^{k});    \label{eq:eg_def}\\
    \vx^{k+1}=&\;\vx^k-\etax g_{\vx}(\vx^{k+1/2};\vy^{k+1/2}),\nonumber\\
    \vy^{k+1}=&\;\vy^k+\etay g_{\vy}(\vx^{k+1/2};\vy^{k+1/2}).\nonumber
\end{align}

From Definition~\ref{def:stackelberg}, we can see that $f$ is approximately a quadratic NC-SC function near $\vz^\ast$. The dynamics of GDA or EG is a linear time invariant system. Specifically, define matrix
\begin{align}
    M=\begin{pmatrix}
    -C& - B\\
    rB^\top & - rA
    \end{pmatrix}.\label{eq:M}
\end{align}
Let $\vz^k=(\vx^k,\vy^k)$. We will see that with an exact gradient oracle, the dynamics of GDA is approximately
\begin{align}
    \vz^{k+1}-\vz^\ast \approx (1+\etax M)\cdot (\vz^k-\vz^\ast).\label{eq:approx_lti_gda}
\end{align}
Similarly, the dynamics of EG is approximately
\begin{align}
    \vz^{k+1}-\vz^\ast\approx (1+\etax M+\etax^2M^2)\cdot (\vz^k-\vz^\ast).\label{eq:approx_lti_eg}
\end{align}

\section{Main results}
\label{sec:local_result}
In this section, we present our main results on the local convergence of GDA and EG.

\subsection{A lower bound of the stepsize ratio for GDA}
\label{sec:lower_bound_ratio}

To start with, we show that for GDA, the stepsize ratio $r$ must be at least $\Omega(\kappa)$ to ensure local convergence to a Stackelberg Equilibrium. Formally, we have the following theorem.
\begin{theorem}
	\label{thm:lower_bound}
	For any given $L,\mu>0$ and $0\le \mux\le L$ such that $\ky=L/\mu\ge 2$, there exists a function $f$ that has a differential Stackelberg Equilibrium with parameters $(L,\mu,\mux)$ such that if choosing a stepsize ratio $r\le \kappa$, the GDA algorithm does not locally converge for all initial points irrespective of any choice of positive stepsize $\etax$.
\end{theorem}

The hard function $f$ in Theorem~\ref{thm:lower_bound} is constructed as $f:\R\times\R\to \R$ with the following expression
\begin{align*}
    {f}(x;y) =& -\frac{\mu}{2}y^2+L xy-\frac{L}{2}  x^2.
\end{align*}
Note that since the above $f$ is a quadratic NC-SC function and global convergence implies local convergence, Theorem~\ref{thm:lower_bound} can be also viewed as a lower bound of the stepsize ratio of GDA for the general NC-SC minimax optimization problem studied in \citep{lin2020gradient}.

For such a simple function $f$, we can directly solve the two eigenvalues of its corresponding $M$ matrix and show that one of them, say $\lambda_1$, has a positive real part. Then $|1+\etax\lambda_1|>1$ for any positive $\etax$ which implies divergence of GDA. However, it does not directly apply to EG because it is much harder to lower bound $|1+\etax\lambda_1+\etax^2\lambda_1^2|$ unless assuming $\etax$ is very small (which is not desirable). We leave the optimal choice of the stepsize ratio for EG as an interesting open problem.

\subsection{Convergence rate of GDA and EG}

Now we show that $r=\Omega(\kappa)$ is sufficient for both GDA and EG to converge. We also provide a tight local convergence rate for them.

\begin{theorem}
\label{thm:main_local_conv}
    Given a function $f$ and its Stackelberg Equilibrium $\vz^\ast$ satisfying all the assumptions in Section~\ref{sec:pre}. Choose $\etax=\frac{1}{4rL}$ and $\etay=\frac{1}{4L}$ with $r\ge 2\kappa$.
    Suppose that $M$ defined in \eqref{eq:M} is diagonalizable and its eigenvalue decomposition is $M=P^{-1}\Lambda P$. Let $C_P=\norm{P}_2\norm{P^{-1}}_2$. There exists some $\delta>0$ such that if $\norm{\vz^0-\vz^\ast}_2\le \delta$, we have the iterates of GDA and EG with an exact gradient oracle satisfy
    \begin{align*}
	\|\vz^T-\vz^\ast\|_2\le C_P\norm{\vz^0-\vz^\ast}_2\left(1-\frac{c_0}{r\kx}\right)^T
	\end{align*}
	for all $T\ge 0$, where $c_0$ is some numerical constant. It means the gradient complexity of achieving an error less than $\epsilon$ is\begin{align*}
	\cO\left(r\kx\log\left(\frac{C_P\norm{\vz^0-\vz^\ast}_2}{\epsilon}\right)\right).
	\end{align*}
\end{theorem}
\paragraph{Remark 1.} The assumption that $M$ is diagonalizable is actually without loss of generality. One can easily obtain similar convergence result without this assumption by using Lemma~\ref{lem:spectral_GDA_EG} and \citep[Chapter 2.1.2,      Theorem~1]{intro_opt}. However, we still keep this assumption because (1) it is a mild assumption as the set of non-diagonalizable matrices has zero measure and (2) it enables us to explicitly derive the dependence on all involved problem-dependent constants like the condition number $C_P$. 

\paragraph{Remark 2.} The factor $2$ in the assumption $r\ge2\ky$ is also \emph{arbitrary}. Following our analysis, we can replace it by $r\ge c\ky$ for any numerical constant $c>1$ and obtain nearly the same convergence result up to constant factors. However, we will always choose $c=2$ in this paper to make expressions as simple as possible.
Even if $c>1$ is not a numerical constant and depends on $\ky$ or $\kx$, local convergence can still be guaranteed although the rate might be worse. Therefore, combined with the negative result in Theorem~\ref{thm:lower_bound}, we know $r=\ky$ is the exact \emph{phase transition point} between convergence and non-convergence.

The above complexity is $\Tilde{\cO}(r\kx)$ which scales linearly in $r$ when $r\ge 2\ky$. So choosing $r=2\ky$ results in the best complexity bound of $\Tilde{\cO}(\ky\kx)$. This rate significantly improves upon the $\Tilde{\cO}(\ky^2\kx)$ complexity under the choice of $r=\Theta(\ky^2)$ as in \citep{lin2020gradient}.

Note that the condition number of the local primal function $\Phi$ is $\kx(\kappa+1)$ as $\Phi$ is $(\kappa+1)L$ smooth. Therefore, the $\Tilde{\cO}(\kappa\kx)$ gradient complexity for $r=2\ky$ is nearly the same as the iteration complexity of running gradient descent directly on $\Phi$ although we do not really have access to its value and gradient!

We also want to point out we fix $\etay=\Theta(1/L)$ and vary $\etax$ for different $r$ because it is well known that $\Theta(1/L)$ is the largest possible step size for $L$ smooth functions to ensure convergence.

To prove Theorem~\ref{thm:main_local_conv}, we note that the dynamics of GDA or EG is approximately a linear time invariant system as in (\ref{eq:approx_lti_gda},\ref{eq:approx_lti_eg}). By a careful analysis of the spectral properties of $M$, we can show that the spectral radius of the transition matrix $I+\etax M$ or $I+\etax M+\etax^2 M^2$ is strictly upper bounded by $1$. In Section~\ref{sec:quad_analysis}, we will introduce our analysis methods in detail through simpler quadratic NC-SC functions. In this case, the dynamical systems~(\ref{eq:approx_lti_gda},\ref{eq:approx_lti_eg}) are exactly LTIs. Then bounding the spectral radius of the transition matrix is sufficient to prove convergence. We will show how to use the results for quadratic functions to prove Theorem~\ref{thm:main_local_conv} in Appendix~\ref{app:main_local_conv}.

It is also worth pointing out that when $f$ is quadratic, it is actually nonconvex-strongly-concave (NC-SC). Also, the results for quadratic NC-SC functions are slightly stronger than Theorem~\ref{thm:main_local_conv}. For example, $M$ no longer needs to satisfy the diaonalizability assumption. As the dynamics of the quadratic functions are globally LTI, a global convergence result is feasible to obtain in this setting. It means that all our local convergence results can be viewed as global results for quadratic NC-SC functions.

Now we provide a matching lower bound to show the convergence rate in Theorem~\ref{thm:main_local_conv} is nearly tight for GDA.
\begin{theorem}
\label{thm:lower_bound_conv_rate_deterministic} Suppose we run GDA or EG with the stepsize choices in Theorem~\ref{thm:main_local_conv} with $r\ge 2\kappa$.
There exists a function $f$ and its Stackelberg Equilibrium $\vz^\ast$ satisfying all the assumptions in Section~\ref{sec:pre} with a diagonalizable $M$ matrix, such that for any there exists a $\vz^0$ that is arbitrarily close to $\vz^\ast$ and the iterates satisfy
\begin{align*}
	\|\vz^T-\vz^\ast\|_2\ge \left(1-\frac{1}{r\kx}\right)^T\norm{\vz^0-\vz^\ast}_2.
	\end{align*}
\end{theorem}

The proof of Theorem~\ref{thm:lower_bound_conv_rate_deterministic} is deferred in Appendix~\ref{app:lower_bound}.
The following theorem establishes the convergence rate of GDA and EG when the gradient oracle is stochastic.

\begin{theorem}
\label{thm:main_local_conv_stochastic}
    With the same assumptions and notation as in Theorem~\ref{thm:main_local_conv}.
     There exists some $\delta>0$ such that if $\norm{\vz^0-\vz^\ast}_2\le \delta$, we have with probability at least $1-p$ for arbitrarily small constant $p$, the iterates of mini-batch SGDA and stochastic EG under Assumption~\ref{ass:zero_mean_bounded_variance} satisfy $\|\vz^T-\vz^\ast\|_2\le \epsilon$ if choosing 
	\begin{align*}
	    T=&O\left(r\kx\log\left(\frac{C_P\norm{\tz^0}_2}{\epsilon}\right)\right)\\
	    S=&O\left(\frac{r^2\kx^2 C_P^2\sigma^2}{L^2\epsilon^2}\log^2\left(\frac{C_P\norm{\tz^0}_2}{\epsilon}\right)\right).
	\end{align*}
	This implies a gradient complexity of 
	\begin{align*}
	    T\cdot S=O\left(\frac{r^3\kx^3 C_P^2\sigma^2}{L^2\epsilon^2}\log^3\left(\frac{C_P\norm{\tz^0}_2}{\epsilon}\right)\right).
	\end{align*}
\end{theorem}

We are not sure whether the above $\tilde{O}\left(\frac{r^3\kx^3}{\epsilon^2}\right)$ complexity is tight or not. Actually, as we will show in next section, the convergence rate for quadratic NC-SC functions in the stochastic setting is $\tilde{O}\left(\frac{r^2\kx^2}{\epsilon^2}\right)$, strictly better than the above theorem. Here we get a worse rate because the deviation of the local dynamics of GDA/EG from the corresponding LTI can depend on the gradient noise, which makes the analysis harder. It is an interesting future direction to see whether the above rate can be improved.

\section{Analysis through quadratic NC-SC functions}
\label{sec:quad_analysis}
In this section, we use the quadratic NC-SC function class as an example to illustrate our analysis.

\subsection{Problem setup for quadratic functions}
We first formulate the problem in the quadratic setting. Consider the quadratic function $f:\R^n\times\R^m\to\R$ of the following form:
\begin{align}
f(\vx;\vy) =& \! - \! \frac{1}{2}(\vy \! - \! \vy^\ast)^\top \! A (\vy \! - \! \vy^\ast) \! + \! (\vx \! - \! \vx^\ast)^\top \! B (\vy \! - \! \vy^\ast)\nonumber \\
& +\frac{1}{2}(\vx-\vx^\ast)^\top C (\vx-\vx^\ast)
\label{eq:quadratic}
\end{align}

for some matrices $A\in\R^{m\times m}$, $B\in\R^{n\times m}$, $C\in\R^{n\times n}$ where $A$ and $C$ are symmetric matrices. $\vz^\ast=(\vx^\ast,\vy^\ast)$ is the differential Stackelberg Equilibrium. In this setting, the assumptions in Section~\ref{sec:pre} are reduced as follows.

\begin{assumption}
	The function $f$ defined in \eqref{eq:quadratic} satisfies
	\begin{enumerate}
		\item $\mu I\le  A\le{L} I$ and $\|B\|_2,\|C\|_2\le {L}$.
		\item $C+BA^{-1}B^\top\ge \mux>0$. 
	\end{enumerate}
\label{assumption:quadratic}
\end{assumption}

Then one can check that the primal function is 
\begin{align*}
\Phi(\vx)=&\frac{1}{2}(\vx-\vx^\ast)^\top \left(C+BA^{-1}B^\top\right)(\vx-\vx^\ast).
\end{align*}
Also, we have $\vx^\ast=\argmin_{\vx} \Phi(\vx)$ and $\vy^\ast=\argmax_{\vy}f(\vx^\ast;\vy)$, i.e., $\vz^\ast$ is a global minimax point for this quadratic NC-SC minimax optimization problem.

\subsection{Convergence in the deterministic setting}
For quadratic functions, (\ref{eq:approx_lti_gda},\ref{eq:approx_lti_eg}) are exactly Linear Time-Invariant (LTI) systems with transition matrices $I+\etax M$ and $I+\etax M+\etax^2M^2$ respectively. Note that from the definition of $M$ in \eqref{eq:M}, these two matrices are non-symmetric and pretty general block matrices. Their eigenvalues could be complex numbers. Although their spectral norms are easy to bound, they could be strictly greater than $1$. Therefore, we have to bound a more fine-grained quantity, the spectral radius, which could be strictly less than the spectral norm. To bound the spectral radius, we provide a careful analysis of the spectral properties of $M$ in the following lemma which is also our main technical contribution.
\begin{lemma}
	Suppose $\lambda=\lambda_0+i\lambda_1$ is an eigenvalue of $M$ where $\lambda_0,\lambda_1\in\R$ and $i=\sqrt{-1}$. Then\begin{enumerate}
		\item $|\lambda_1| \le \sqrt{r}{L}$.
		\item If $\lambda_1\neq 0$, we have $\lambda_0\le-\frac{\mu(r-\ky)}{2}$. 
		\item If $r\ge1$, the spectral norm of $M$ is bounded: $\lambda_0^2+\lambda_1^2\le \norm{M}_2^2\le 4r^2{L}^2$.
		\item Choosing $r>\kappa$, if $\lambda_1=0$, we have $\lambda_0\le -\mux$.
		\item If $r> \kappa$ and $\mux>0$, we must have $\lambda_0 < 0$.
	\end{enumerate} 
	\label{lem:spectral_M}
\end{lemma}
Remark that Lemma~\ref{lem:spectral_M}.1 and~\ref{lem:spectral_M}.5 together suggest that if $\mux>0$ and $r>\ky$, every eigenvalue $\lambda$ of $M$ has a strictly negative real part and a bounded imaginary part. Then if fixing $r=\etay/\etax>\ky$ and choosing a small enough $\etax$ (or $\etay$), one can guarantee $|1+\etax \lambda|<1$, i.e., the spectral radius of the transition matrix of GDA is less than $1$ which implies convergence. Similary, the convergence of EG can also be guaranteed.

Lemma~\ref{lem:spectral_M} also suggests that the angle between the complex number $\lambda$ and the real line in the complex plane is bounded by some $\theta<\pi/2$. If further assuming $r\ge c\ky$ for some numerical constant $c>1$, we can get a tight bound of $\theta$ which allows us to choose $\etax$ to make $|1+\etax \lambda|$ as small as possible. Then we can obtain a quantitative and tight bound of the spectral radius of the transition matrices in the following lemma.

\begin{lemma}
	\label{lem:spectral_GDA_EG}
	Choose $\etax=\frac{1}{4r{L}}$ and $\etay=\frac{1}{4{L}}$ where $r\ge2\kappa$. Let $\rho_1$ be the spectral radius of $I+\etax M$ and $\rho_2$ be that of $I+\etax M+\etax^2 M^2$. We have
	\begin{align*}
	\max\{\rho_1,\rho_2\}\le 1-\frac{1}{64r\kx}.
	\end{align*}
\end{lemma}

\begin{figure*}[t]
     \centering
     \begin{subfigure}[b]{0.49\columnwidth}
         \centering
         \includegraphics[width=0.9\textwidth]{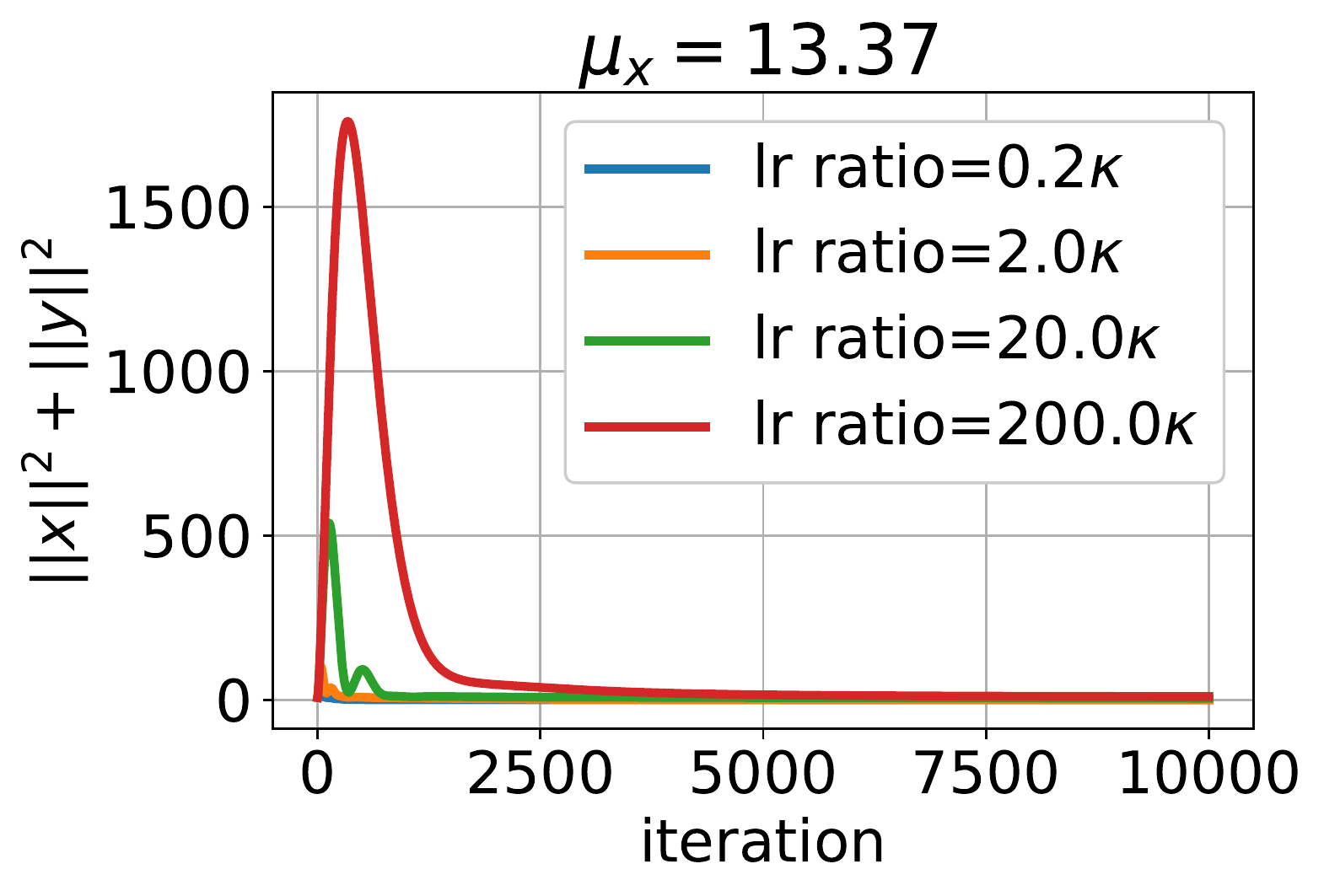}
         \caption{}
         \label{fig:quad_training_a}
     \end{subfigure}
     \hfill
     \begin{subfigure}[b]{0.49\columnwidth}
         \centering
         \includegraphics[width=0.9\textwidth]{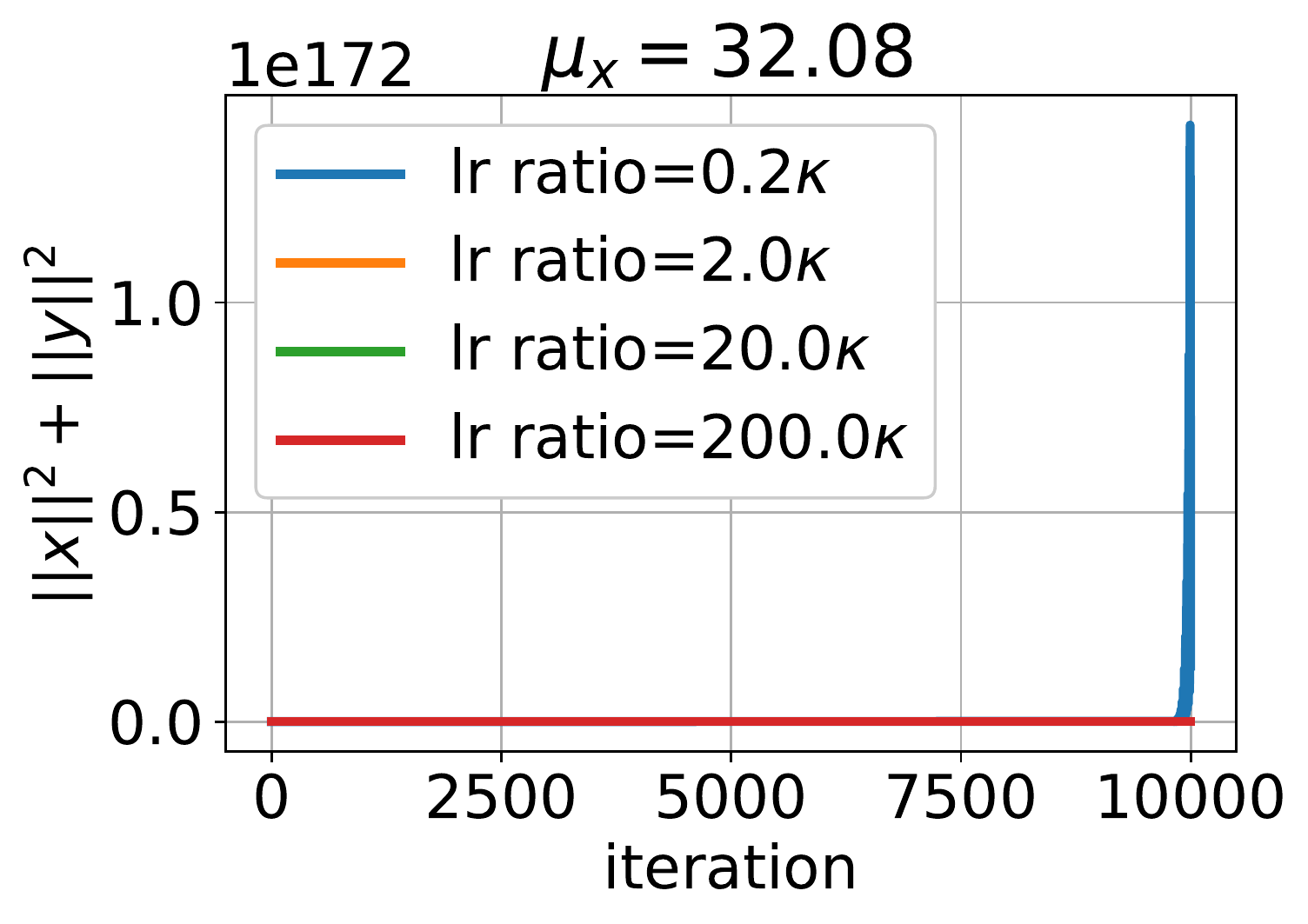}
         \caption{}
         \label{fig:quad_training_b}
     \end{subfigure}
        \caption{Training curves of GDA on two independently generated quadratic functions satisfying Assumption~\ref{assumption:quadratic} with four different stepsize ratios. The corresponding $\mux$ value is shown in the title of each plot.}
        \label{fig:quad_training}
\end{figure*}

\begin{figure*}[t]
     \centering
     \begin{subfigure}[b]{0.48\textwidth}
         \centering
         \includegraphics[width=0.9\textwidth]{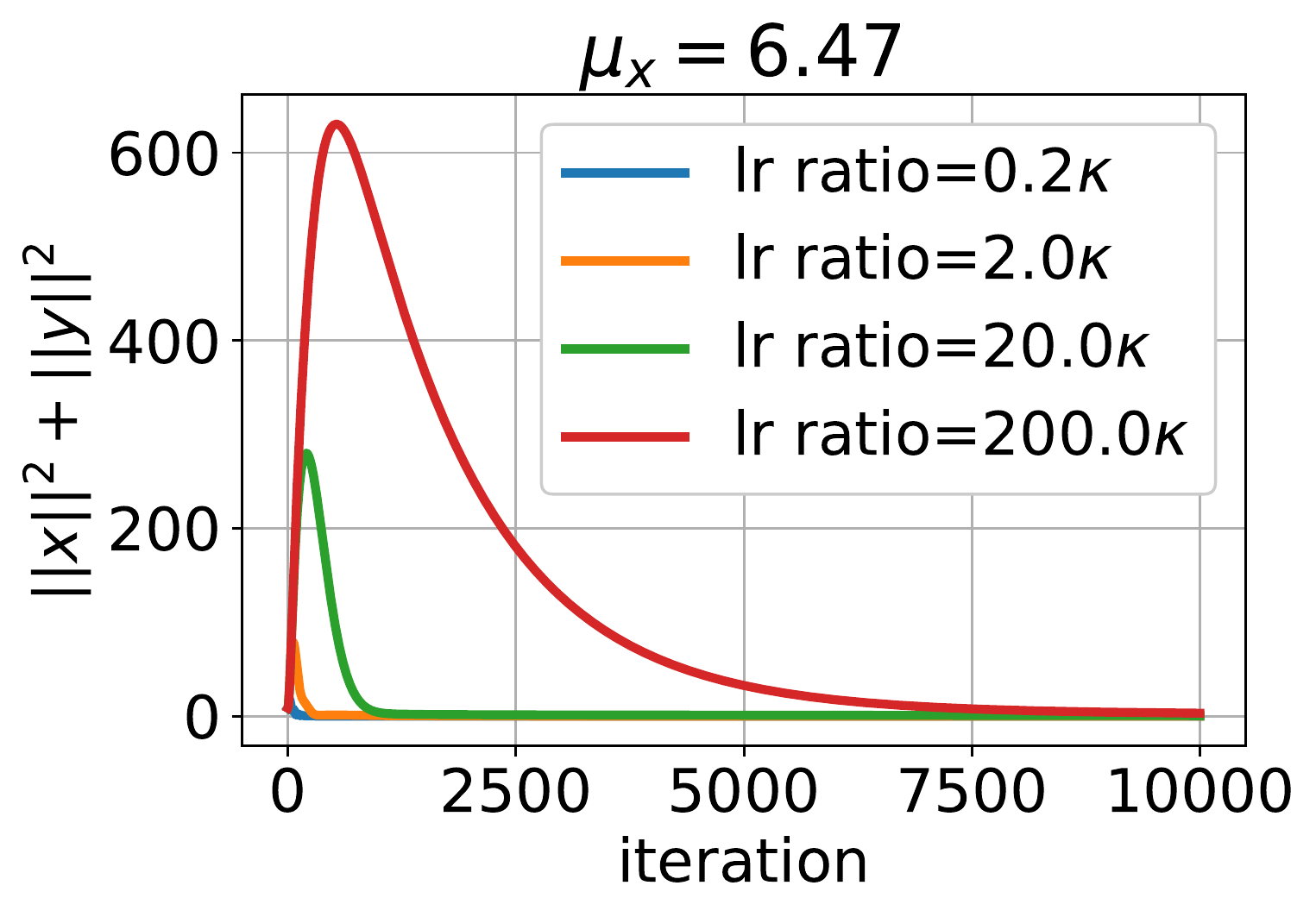}
         \caption{}
         \label{fig:nonquad_training_a}
     \end{subfigure}
     \hfill
     \begin{subfigure}[b]{0.48\textwidth}
         \centering
         \includegraphics[width=0.9\textwidth]{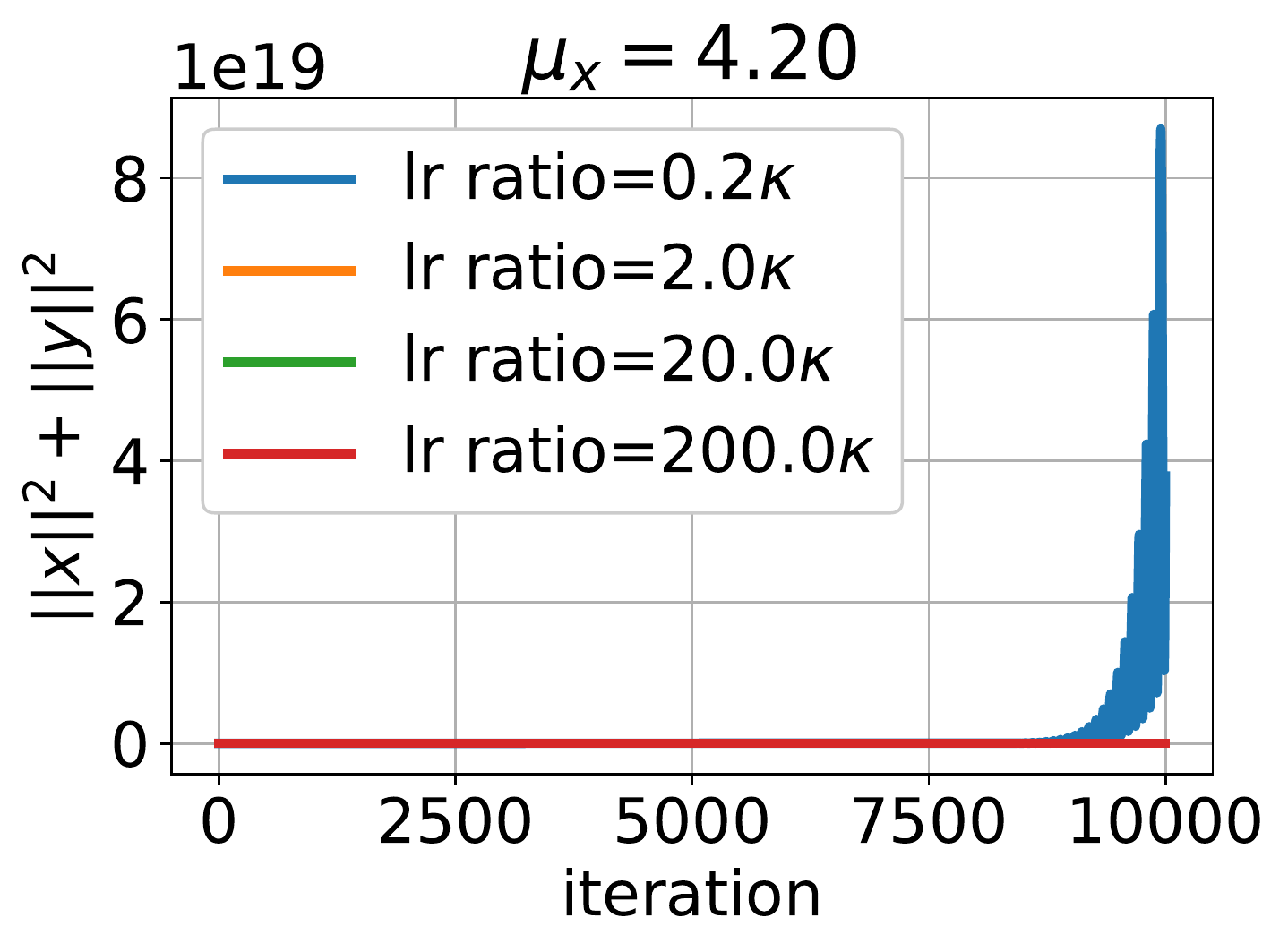}
         \caption{}
         \label{fig:nonquad_training_b}
     \end{subfigure}
        \caption{Training curves of GDA on two independently generated non-quadratic functions defined in \eqref{eq:nonquad_example} with four different stepsize ratios. The corresponding $\mux$ value is shown in the title of each plot.}
        \label{fig:nonquad_training}
\end{figure*}

\begin{figure*}[t]
     \centering
     \begin{subfigure}[b]{0.49\textwidth}
         \centering
         \includegraphics[width=0.9\textwidth]{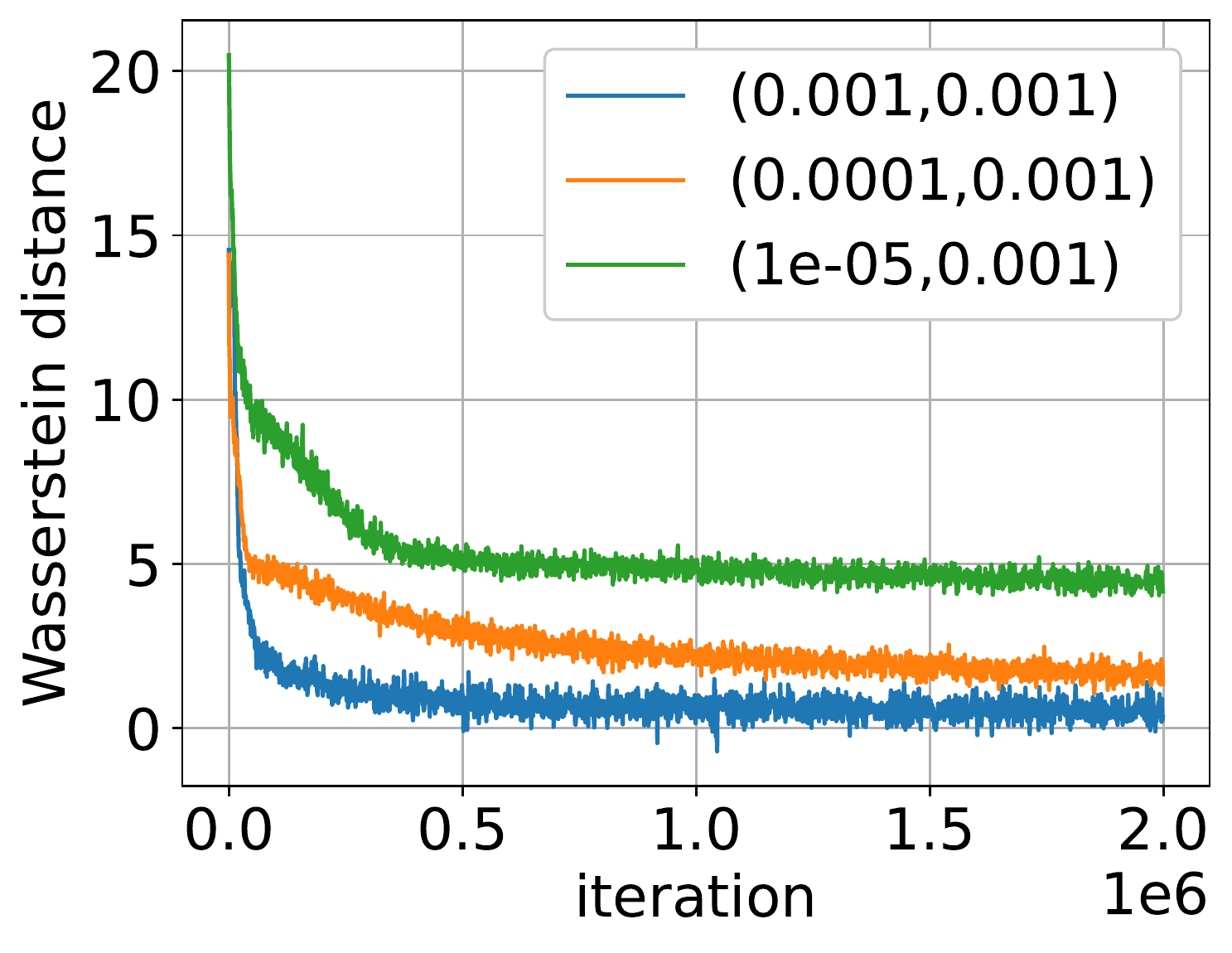}
         \caption{}
         \label{fig:mot_exp_a}
     \end{subfigure}
     \hfill
     \begin{subfigure}[b]{0.49\textwidth}
         \centering
         \includegraphics[width=0.9\textwidth]{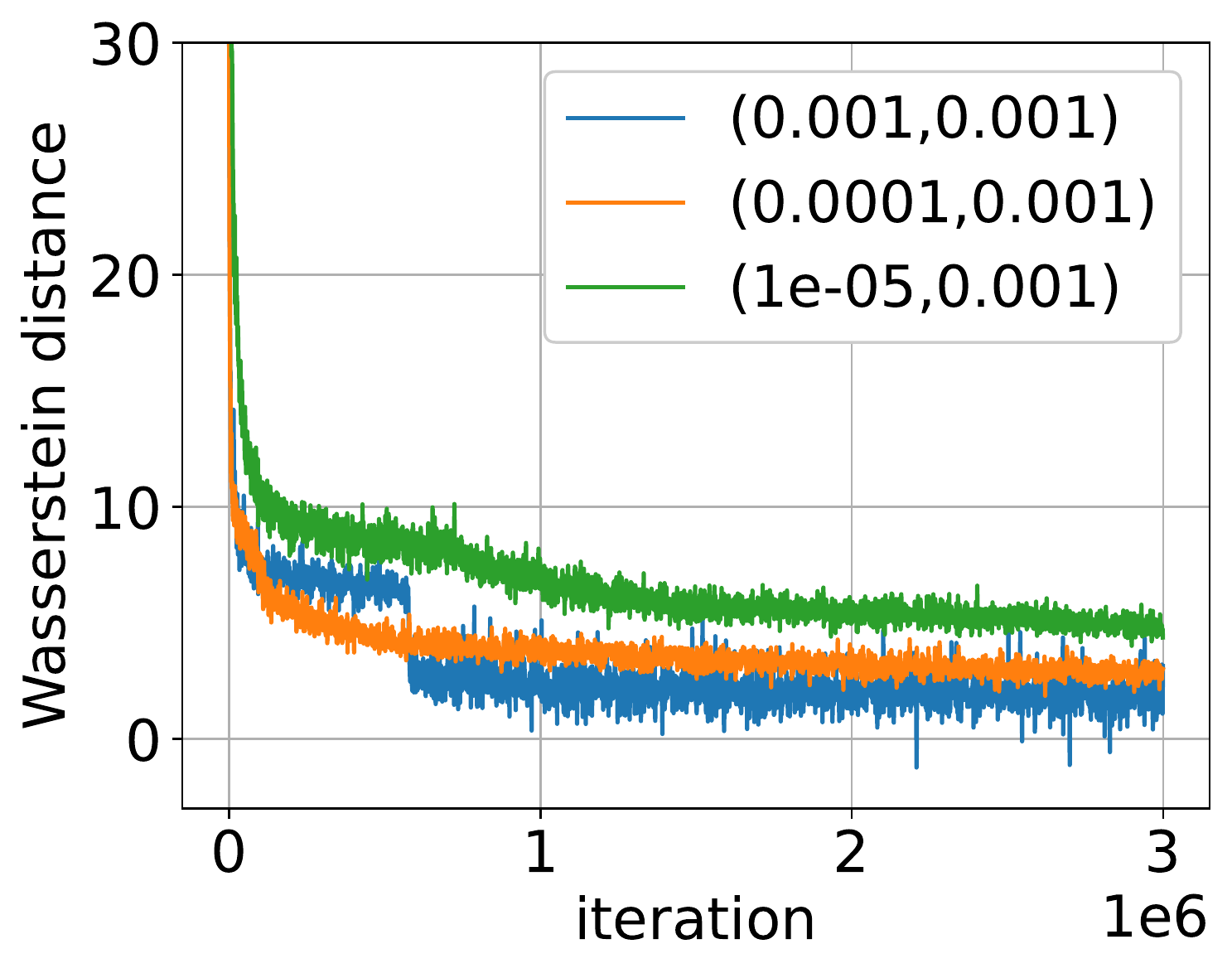}
         \caption{}
         \label{fig:mot_exp_b}
     \end{subfigure}
        \caption{Evolution of Wasserstein Distance on MNIST (a) and CIFAR10 (b). For both MNIST and CIFAR10, we train WGAN-GP models using under different stepsize raios $(\etax,\etay)=(0.001,0.001),(0.0001,0.001),(1e-05,0.001)$.}
        \label{fig:mot_exp}
\end{figure*}

The proofs of Lemma~\ref{lem:spectral_M} and Lemma~\ref{lem:spectral_GDA_EG} are deferred in Appendix~\ref{app:spectral}.
With Lemma~\ref{lem:spectral_GDA_EG}, 
the convergence result of GDA and EG for quadratic NC-SC functions follows naturally.

\begin{theorem}
	Let $M=P^{-1}J P$ be its Jordan decomposition. Define $C_P=\|P\|_2\|P^{-1}\|_2$ as the condition number of $P$ and $s$ be the size of the largest Jordan block in $J$.
	Under Assumption~\ref{assumption:quadratic}, choosing $\etax=\frac{1}{4rL}$ and $\etay=\frac{1}{4L}$ with $r\ge 2\ky$, we have the iterates of GDA or EG with an exact gradient oracle satisfy
	\begin{align*}
	\|\vz^T \! - \! \vz^\ast\|_2 \! \le \!  sC_P \! \norm{\vz^0 \! - \! \vz^\ast}_2 \! T^{s-1} \! \left(1 \! - \! \frac{c_0}{r\kx}\right)^{T-s+1}
	\end{align*}
	 for all $T\ge0$, where $c_0$ is some numerical constant. It means the gradient complexity of achieving $\|\vz^T-\vz^\ast\|_2\le\epsilon$ is\begin{align*}
\cO\left(sr\kx \log\left(\frac{sr\kx C_P\norm{\vz^0-\vz^\ast}_2}{\epsilon}\right)\right).
	\end{align*}
	\label{thm:converence_nond_1}
\end{theorem}

Note that Theorem~\ref{thm:converence_nond_1} is a global convergence result and does not assume $M$ is diagonalizable.
However, the gradient complexity has a factor of $s$ which could be proportional to the dimension $d=m+n$ in the worst case when $M$ is not diagonalizable but is usually very small.

\subsection{Convergence in the stochastic setting}
The following theorem establishes a convergence result for mini-batch SGDA. 
\begin{theorem}
	\label{thm:conv_general_sgda}
	We use the same assumptions and parameters as in Theorem~\ref{thm:converence_nond_1} except that now we consider stochastic gradient oracles. The iterates of SGDA and stochastic EG satisfy  
	\begin{small}
	\begin{align*}
	\E[\|\vz^T\!-\!\vz^\ast\|_2^2]\le   C_P^2\left(1\!-\!\frac{c_1}{r\kx}\right)^{2T}\!\norm{\vz^0-\vz^\ast}_2^2\!+\!\frac{c_2r\kx C_P^2\sigma^2}{L^2S}
	\end{align*}
	\end{small}
	 for all $T\ge0$, where $c_1$ and $c_2$ are numerical constants. It means the gradient complexity of achieving $\E[\|{\vz^T-\vz^\ast}\|_2]\le\epsilon$ is\begin{align*}
	\cO\left(r\kx\log\left(\frac{C_P\norm{\vz^0-\vz^\ast}_2}{\epsilon}\right)\max\left\{1,\frac{r\kx C_P^2\sigma^2}{L^2 \epsilon^2}\right\}\right).
	\end{align*}
\end{theorem}

The proofs of Theorem~\ref{thm:converence_nond_1} and Theorem~\ref{thm:conv_general_sgda} can be found in Appendix~\ref{app:quad_conv}.
As we mentioned before, the above $\Tilde{O}(\frac{r^2\kx^2}{\epsilon^2})$ rate is better than the stochastic local convergence counterpart in Theorem~\ref{thm:main_local_conv_stochastic}.

Again note that when $r=2\kappa$, the rate is $\tilde{O}\left(\frac{(\ky\kx)^2}{\epsilon^2}\right)$ which is the gradient complexity of stochastic gradient descent directly on $\Phi$ although we do not have access to the stochasic oracle of $\Phi$!

\section{Experimental results}
\label{sec:exp}

\subsection{Quadratic functions}
\label{sec:exp_quad}
In this section, we empirically show the convergence behavior of GDA with different stepsize ratios on quadratic functions. We set $\vz^\ast=0$ w.l.o.g. as it does not affect the convergence behavior. The matrices $A,B,C\in\R^{4\times 4}$ are randomly generated and processed to satisfy Assumption~\ref{assumption:quadratic}. We choose $L=100$, $\mu=1$ in the beginning and compute $\mux$ after the matrices are sampled. We keep $\etay=\frac{1}{2L}$ and change the stepsize ratio by varying $\etax$.

Figure~\ref{fig:quad_training} shows the training curves under two independent samples for four different stepsize ratios. The orange curve uses the best stepsize ratio, $2\ky$, as suggest by our theorems. The red curve uses the stepsizes as in \citep{lin2020gradient} (note $200\kappa=2\kappa^2$). According to Figure~\ref{fig:quad_training_a}, as the stepsize ratio increases from $2\ky$, the convergence becomes slower. For the red curve, since $\etax$ is much smaller than $\etay$, $\vx$ updates much more slowly than $\vy$. As a result, the training can be roughly splitted into two phases. During the first phase (first $1000$ iterations in Figure~\ref{fig:quad_training_a}), the change of $\vx$ can be ignored and $\vy$ is updated to maximize $f$. Therefore the red curve goes up. During the second phase (after $1000$ iterations in Figure~\ref{fig:quad_training_a}), $\vy$ is updated to keep maximizing $f$ for each $\vx$. Therefore essentially, $\vx$ is trained to minimze the primal function $\Phi(\vx)$ and thus the curve goes down. However, such phenomenon does not exist for the orange curve, as both variables are basically simultaneously trained to the optimal point.

For the blue curve which uses a stepsize ratio even much smaller than $2\ky$, although it converges even faster than the orange curve as in Figure~\ref{fig:quad_training_a}, it sometimes does diverge as in Figure~\ref{fig:quad_training_b}. This is consistent with our lower bound in Theorem~\ref{thm:lower_bound}.

\subsection{A non-quadratic function}
\label{sec:exp_nonquad}
Although our theory only applies to quadratic functions or local convergence, the empirical training behavior we observed above may exist for more general functions. In this section, we empirically run GDA on a non-quadratic function defined as follows. We first generate an instance of quadratic funtions, $f_0$, as in Section~\ref{sec:exp_quad}. Then we add a non-quadratic function of $\vx$ only. Formally, we run GDA on the following function:
\begin{align}
    f(\vx;\vy):=&f_0(\vx;\vy)+\tfrac{L}{n}\textstyle\sum_{i=1}^n\left[\log\left(1+\exp(a(x_i-b_i)\right)\right.\nonumber\\&\left.+\log\left(1+\exp(-a(x_i-b_i)\right)\right].
    \label{eq:nonquad_example}
\end{align}
The training curves for different stepsize ratios are shown in Figure~\ref{fig:nonquad_training}, with similar behaviors to those in Figure~\ref{fig:quad_training}.

\subsection{GAN training}
\label{sec:mot_exp}

In this section, we empirically show the convergence behavior of GANs trained using GDA with different stepsize ratios. The WGAN-GP model in \citep{Gulrajani2017ImprovedTO} was trained using ADAM, a variant of GDA, with $\etax=\etay=0.0001$. Moreover, they ran multiple steps of inner maximization each outer minimization step. Our experiment uses the same WGAN-GP model. However, different from their algorithm, we use simultaneous GDA with the same number of gradient steps for both variables. Figure~\ref{fig:mot_exp} shows the generated images of the learned generator and the evolution of the Wasserstein Distance which is also the primal function $\Phi$. As we can see, simultaneous GDA with $\etax=\etay=0.001$ is able to converge with a high speed for both MNIST and CIFAR10. However, if we decrease $\etax$, i.e., increase the stepsize ratio, it converges with a lower speed. The convergent value is also larger. 
We will show more experimental results for different models under different stepsize ratios in the Appendix.

\section{Conclusion}
\label{sec:conclusion}
There exists a wide gap between existing theory on gradient-based minimax optimization which suggests a significantly large max/min stepsize ratio and state-of-the-art GAN training algorithms where similar learning rates are used for the min and max players. In this paper, we take one step towards understanding and bridging this gap by providing a tight local convergence analysis of GDA near a strict Stackelberg Equilibrium. In particular, we characterize the optimal choices of the stepsize values and resulting stepsize ratio and further provide a local convergence guarantee under these stepsize selection. We also extend the analysis to extra-gradient and stochastic GDA methods. 

There remain several open problems and future directions following our paper. First, in the stochastic setting, the complexity bounds of GDA and EG for local convergence in Theorem~\ref{thm:main_local_conv_stochastic} are larger than the bounds under quadratic objective functions in Theorem~\ref{thm:conv_general_sgda}. Hence, an interesting future direction is to improve the bounds in Theorem~\ref{thm:main_local_conv_stochastic}. In addition, the extension of our results to optimistic GDA (OGDA) methods will be another interesting problem. 
Finally, our analysis focuses on local convergence and quadratic NC-SC objective functions. The problem of characterizing the best stepsize ratio for general NC-SC functions and NC-NC functions still remains open.

\section{Acknowledgement}
This work was supported, in part, by the MIT-IBM Watson AI Lab, ONR Grant N00014-20-1-2394, a Vannevar Bush fellowship from Office of the Secretary of Defense, and a CUHK Direct Grant for Research. The authors thank the anonymous reviewers for their constructive feedback and suggestions.

\bibliographystyle{plainnat}
\bibliography{arxiv_version}

\newpage

%
%




%

%

\onecolumn

\section{Proofs for lower bounds}
\label{app:lower_bound}
In this section, we provide the proofs for all the lower bounds. Specifically, we prove Theorem~\ref{thm:lower_bound} and Theorem~\ref{thm:lower_bound_conv_rate_deterministic}.
\subsection{Proof of Theorem~\ref{thm:lower_bound}}
\begin{proof}[Proof of Theorem~\ref{thm:lower_bound}]
	Consider the following quadratic NC-SC function ${f}:\R\times\R\to \R$:
	\begin{align*}
		{f}(x;y) =& -\frac{\mu}{2}(y-\kappa x)^2 + \frac{\mu}{2}(\ky^2-\ky)x^2\\
		=&-\frac{\mu}{2}y^2+L xy-\frac{L}{2}  x^2\\
		=&\frac{\mu}{2} \begin{pmatrix}
		x & y
		\end{pmatrix}\begin{pmatrix}
		-\ky & \ky\\
		\ky & -\mu
		\end{pmatrix} \begin{pmatrix}
		x \\ y
		\end{pmatrix}.
	\end{align*}
	Note that ${\Phi}(x)=\max_{y\in\R} {f}(x;y)=\frac{\mu}{2}(\ky^2-\ky)x^2\ge\frac{\mux}{2}x^2$. It is easy to verify that $f$ satisfies Assumption~\ref{assumption:quadratic} and $(0,0)$ is its unique differential Stackelberg Equilibrium with parameters $(L,\mu,\mux)$. Indeed, it is also the global minimax point. The $M$ matrix for $f$ is 
	\begin{align*}
	    M=\begin{pmatrix}
	\ky & -\ky\\
	r \ky & -r
	\end{pmatrix}.
	\end{align*}
	Let $\lambda_1$ and $\lambda_2$ be the eigenvalues of $M$, we have
	\begin{align*}
	\tr(M)=\lambda_1+\lambda_2 = \ky-r.
	\end{align*}
	There are two possible cases:
	\begin{itemize}
	    \item If $r<\kappa$, we know $\lambda_1+\lambda_2>0$ which implies at least one of them has a positive real part. Suppose $\re[\lambda_1]>0$ without loss of generality. We know $|1+\etax\lambda_1|>1$ for every $\etax>0$.
	    \item If $r=\ky$, $M=\begin{pmatrix}
	\ky & -\ky\\
	 \ky^2 & -\ky
	\end{pmatrix}$. We can directly compute its eigenvalues $\lambda_1=i\sqrt{\ky^3-\ky^2}$ and $\lambda_2=-i\sqrt{\ky^3-\ky^2}$ where $i=\sqrt{-1}$. We also have $|1+\etax\lambda_1|=\sqrt{1+\etax^2(\ky^3-\ky^2)}>1$ for every $\etax>0$.
	\end{itemize}
	  Then the power sequence $\{(I+\etax M)^k\}_{k=0}^{\infty}$ will diverge and we complete the proof.
\end{proof}

\subsection{Proof of Theorem~\ref{thm:lower_bound_conv_rate_deterministic}}
\begin{proof}[Proof of Theorem~\ref{thm:lower_bound_conv_rate_deterministic}]
Consider the following quadratic NC-SC function $f:\mathbb{R}\times\mathbb{R}\to \mathbb{R}$
\begin{align*}
	f(x,y):=-\frac{1}{2}L x^2+bxy-\frac{1}{2}\mu y^2,
	\end{align*}
where $b = \sqrt{\mu (L+\mu_x)}$. The primal function is $\Phi(x):=\max_y f(x,y)=\tfrac{1}{2}\mu_x x^2$. Therefore, it is straight-forward to verfity that $f$ satisfies Assumption~\ref{assumption:quadratic} and $(0,0)$ is its unique differential Stackelberg Equilibrium with parameters $(L,\mu,\mux)$. Assume $\kappa=L/\mu\ge 2$. We first show the lower bound of GDA to achieve $(0,0)$. The $M$ matrix is 
\begin{align*}
M=\begin{pmatrix}
L & -b\\ r b & -\mu r 
\end{pmatrix}.
\end{align*}
 Let $\lambda_1$ and $\lambda_2$ be the two eigenvalues of $M$, we have
\begin{align*}
	\lambda_{1,2} = -\frac{1}{2}\left(\mu r-L\right) \pm \frac{1}{2}\sqrt{\left(\mu r-L\right)^2-4r\mu\mu_x}.
\end{align*}

\begin{align*}
0\ge	\lambda_1 &= -\frac{1}{2}(\mu r-L)\left(1-\sqrt{1-\frac{4r\mu\mu_x}{(\mu r-L)^2}}\right)\\
& \ge-\frac{2r\mu\mux}{r\mu-L}\ge -4\mux.
\end{align*}

Let $s_1$ be the corresponding eigenvalue of $I+\eta_x M$, it satisfies
\begin{align*}
	 0\le 1-\frac{\mu_x}{rL}\le s_1=1+\eta_x \lambda_1\le 1.
\end{align*}

We adversarially choose the initial point $\vz^0$ such that it is parallel to the eigenvector of $I+\eta_x M$ corresponding to $s_1$. Note that for any $\delta>0$, we can choose such a $\vz^0$ such that $\norm{\vz^0-\vz^\ast}_2\le \delta$. Then we have
\begin{align*}
 \vz^{T} = (I+\etax M)^T\vz^0 = s_1^T\vz^0,
\end{align*}
which implies
\begin{align*}
    \norm{\vz^T}_2=s_1^T\norm{\vz^0}_2\ge \left(1-\frac{1}{r\kx}\right)^T\norm{\vz^0}_2.
\end{align*}
For EG, note that $1+\etax\lambda_1^2+\etax^2\lambda_2^2\ge 1+\etax\lambda_1^2=s_1$. The rest of proof is the same as that of GDA.
\end{proof}

\section{Analysis of the spectral properties}
\label{app:spectral}

In this section, we prove the spectral properties of $M$ (Lemma~\ref{lem:spectral_M}) and the spectrual radius bound of GDA and GD (Lemma~\ref{lem:spectral_GDA_EG}).
We first present the following lemma and its corollary useful for analyzing spectral properties of a matrix.
\begin{lemma}
	Suppose $G=A+iB\in \C^{d\times d}$ where $A,B\in \R^{d\times d}$ are both symmetric real matrices and $i=\sqrt{-1}$. If there exist $a,b\in \R$ such that $aA+bB$ is strictly positive or negative definite, then $G$ is invertible, i.e., $\det(G)\neq 0$.
	\label{lem:1}
\end{lemma}

\begin{proof}[Proof of Lemma~\ref{lem:1}]
	If $\det(G)=0$, we know there exists $z=z_0+iz_1\neq 0$ for some $z_0,z_1\in\R^d$ such that $Gz=0$. Then we have\begin{align*}
	z^\dag G z = z_0^\top A z_0+z_1^\top A z_1+i\left(z_0^\top B z_0+z_1^\top B z_1\right)=0,
	\end{align*}
	where $z^\dag$ is the conjugate transpose of $z$. It means both the real and imaginary parts of $z^\dag G z$ are zero. So is any linear combination of them. Therefore\begin{align*}
	z_0^\top (aA+bB) z_0+z_1^\top (aA+bB) z_1 = 0
	\end{align*}
	which contradicts with the strict definiteness of $aA+bB$. Therefore $\det(G)\neq 0$.
\end{proof}
Choosing $a=1,b=0$ or $a=0,b=1$, we obtain the following corollary immediately.
\begin{corollary}
	Suppose $G=A+iB\in \C^{d\times d}$ where $A,B\in \R^{d\times d}$ are both symmetric real matrices and $i=\sqrt{-1}$. If at least one of $A$ and $B$ is strictly positive or negative definite, then $G$ is invertible, i.e., $\det(G)\neq 0$.
	\label{cor:1}
\end{corollary}

Now we are ready to prove Lemma~\ref{lem:spectral_M}, the spectral properties of $M$.
\begin{proof}[Proof of Lemma~\ref{lem:spectral_M}]
Before proving the five parts separately, we will first prove that $\lambda_0\le 0$ if $r>\ky$, which is a slightly weaker version of part 5 (where we claim $\lambda_0< 0$ if $r>\ky$). We will prove it by contradiction.
	 Assume that $r>\ky$ but $\lambda_0>0$. The eigenvalue satisfies $\det(\lambda I-M)=0$. Note by the theory of Schur Complement,
		\begin{align*}
		\det(\lambda I-M)&=\det \begin{pmatrix}
		C+\lambda I & B\\
		-r B^\top & r A+\lambda I
		\end{pmatrix}\\
		&= \det\left(r A+\lambda I \right)\det\left( C+\lambda I+B \left(A+\frac{\lambda}{r}I\right)^{-1}B^\top \right),
		\end{align*}
		where in the second equality we use the fact that $r A+\lambda I$ is invertible due to Corollary~\ref{cor:1}. Let $A=U\Lambda U^\top$ where $U\in\R^{m\times m}$ is orthonormal and $\Lambda=\text{diag}(\sigma_1,\ldots,\sigma_m)$ is diagonal.  We know $\mu\le \sigma_i \le {L}$ for every $i\in[m]$. Let $\widetilde{B}=BU$. Define
		\begin{align*}
		H(\lambda) :=& C+\lambda I+B \left(A+\frac{\lambda}{r}I\right)^{-1}B^\top\\
		=& C+\lambda I+\widetilde{B}D\widetilde{B}^\top
		\end{align*}
		where $D=\text{diag}(d_1,\ldots,d_m)$ is a diagonal matrix with\begin{align*}
		d_i = \frac{ \sigma_i+\frac{\lambda_0}{r}-i\left(\frac{\lambda_1}{r}\right)
		}{\left(\sigma_i+\frac{\lambda_0}{r}\right)^2+\left(\frac{\lambda_1}{r}\right)^2}
		\end{align*}whose real part is strictly positive. We know $\det(\lambda I-M)=0$ if and only if $\det(H(\lambda))=0$ as $\det(rA+\lambda I)\neq 0$.
		The real part of $H(\lambda)$ is\begin{align*}
		\re[H(\lambda)]= C+\lambda_0 I+ \widetilde{B}\re[D]\widetilde{B}^\top.
		\end{align*}
		Since $C\ge -{L}I$ and $\re[D]\ge 0$, we know that if $\lambda_0>{L}$, then $\re[H(\lambda)]>0$. Then by Corollary~\ref{cor:1}, $\det(H(\lambda))\neq 0$. Therefore we must have $\lambda_0\le {L}$. 
		
		Now suppose $\lambda_1\neq 0$. Then note the imaginary part of $H(\lambda)$ is\begin{align*}
		\im[H(\lambda)]=\lambda_1 \left(I-\frac{1}{r}\widetilde{B}E\widetilde{B}^\top\right),
		\end{align*}
		where $E=\text{diag}(e_1,\ldots,e_m)$ with \begin{align*}
		e_i=\frac{ 1
		}{\left(\sigma_i+\frac{\lambda_0}{r}\right)^2+\left(\frac{\lambda_1}{r}\right)^2}.
		\end{align*}
		Note\begin{align*}
		\re[H(\lambda)]+\frac{r\mu}{\lambda_1}\im[H(\lambda)]=C+\lambda_0I+r\mu I+\widetilde{B}F\widetilde{B}^\top,
		\end{align*}	
		where $F=\text{diag}(f_1,\ldots,f_m)$ with \begin{align*}
		f_i=\frac{ \sigma_i+\frac{\lambda_0}{r}-\mu
		}{\left(\sigma_i+\frac{\lambda_0}{r}\right)^2+\left(\frac{\lambda_1}{r}\right)^2}>0.
		\end{align*}
		Then we have $\re[H(\lambda)]+\frac{r\mu}{\lambda_1}\im[H(\lambda)]>0$ which contradicts with Lemma~\ref{lem:1}. Therefore we must have $\lambda_1=0$, i.e., $\lambda$ is real.
		
		Then we can calculate
		\begin{align*}
		H(\lambda) &= C+\lambda_0 I+B \left(A+\frac{\lambda_0}{r}I\right)^{-1}B^\top \\
		&=C+\lambda_0 I+B \left(A+\frac{\lambda_0}{r}I\right)^{-1}\left[\left(A+\frac{\lambda_0}{r}I\right)-\frac{\lambda_0}{r}I\right]A^{-1}B^\top \\
		&=C+\lambda_0 I+BA^{-1}B^\top -\frac{\lambda_0}{r}B \left(A+\frac{\lambda_0}{r}I\right)^{-1} A^{-1}B^\top\\
		&\ge C+BA^{-1}B^\top+\lambda_0 I-\frac{\lambda_0}{r\mu}BA^{-1}B^\top\\
		&=\left(1-\frac{\lambda_0}{r\mu}\right)\left( C+BA^{-1}B^\top\right)+\lambda_0\left( I+\frac{C}{r\mu}\right)\\
		&>0,
		\end{align*}
		which again contradicts with Corollary~\ref{cor:1}, where in the last inequality we use the fact $1-\lambda_0/r\mu>0$ and $I+C/r\mu>0$ if $r>\kappa$ and $0<\lambda_0\le L$. Therefore we must have $\lambda_0\le 0$.
		
		Now we are ready to prove the five parts one by one.
	\begin{enumerate}
		\item 
		First note that if $\abs{\lambda_1}>\sqrt{r}{L}$, we know $r A+\lambda I$ is still invertible without assuming $\lambda_0>0$ as $\lambda_1\neq 0$, which means we are still able to define $H(\lambda)$ and obtain \begin{align*}
		\im[H(\lambda)]=\lambda_1 \left(I-\frac{1}{r}\widetilde{B}E\widetilde{B}^\top\right).
		\end{align*}
		It is easy to see $\|E\|_2\le \left(\frac{r}{\lambda_1}\right)^2$.
		Noting $\|\widetilde{B}\|_2=\|B\|_2\le{L}$, we have\begin{align*}
		I-\frac{1}{r}\widetilde{B}E\widetilde{B}^\top>0.
		\end{align*}
		Then by Corollary~\ref{cor:1}, $\det(H(\lambda))\neq 0$ which leads to contradiction. Therefore we must have $|\lambda_1|\le \sqrt{r}{L}$.
		
		\item 
		Note if $\lambda_1\neq 0$,  we have
		\begin{align*}
		\re[H(\lambda)]+\frac{\alpha r\mu}{\lambda_1}\im[H(\lambda)]=C+\lambda_0I+\alpha\mu r I+ \widetilde{B}F^{(\alpha)} \widetilde{B}^\top,
		\end{align*}
			where we choose $\alpha=\frac{r+\ky}{2r}$ and $F^{(\alpha)}=\text{diag}(f^{(\alpha)}_1,\ldots,f^{(\alpha)}_m)$ with \begin{align*}
		f^{(\alpha)}_i=\frac{ \sigma_i+\frac{\lambda_0}{r}-\alpha\mu
		}{\left(\sigma_i+\frac{\lambda_0}{r}\right)^2+\left(\frac{\lambda_1}{r}\right)^2}.
		\end{align*}
		Suppose $\lambda_0>-(\alpha \mu r-L)=-\frac{\mu(r-\ky)}{2}$. We have $C+\lambda_0I+\alpha\mu r I>0$.
	Also note that
		\begin{align*}
		    \sigma_i+\frac{\lambda_0}{r}-\alpha\mu\ge \sigma_i-  2\alpha\mu+L/r\ge (1-2\alpha+\kappa/r)\mu=0
		\end{align*}
		which implies $f^{(\alpha)}_i\ge 0$ and thus $F^{(\alpha)}$ is positive semi-definite. Therefore $	\re[H(\lambda)]+\frac{\alpha r\mu}{\lambda_1}\im[H(\lambda)]>0$ which leads to contradiction by Lemma~\ref{lem:1}. Therefore $\lambda_0\le-\frac{\mu(r-\ky)}{2}$.
		\item 
		Note that the spectral radius of a matrix is always bounded by its operator norm, i.e., $\lambda_0^2+\lambda_1^2\le \|M\|^2_2$. Note\begin{align*}
		\left\|M \begin{pmatrix}
		x \\ y
		\end{pmatrix}\right\|_2 &=\left\| \begin{pmatrix}
		-Cx-By \\ r B^\top x-r Ay
		\end{pmatrix} \right\|_2\\
		&= \sqrt{\|-Cx-By\|_2^2+\|r B^\top x-r Ay\|_2^2}\\
		&\le \sqrt{2\left(\|Cx\|_2^2+\|By\|_2^2  \right)+2r^2\left(\| B^\top x\|_2^2+\| Ay\|_2^2\right)}\\
		&\le \sqrt{2(r^2+1){L}^2 \left(\|  x\|_2^2+\| y\|_2^2\right) }\\
		&\le 2r{L} \left\| \begin{pmatrix}
		x \\ y
		\end{pmatrix}\right\|_2.
		\end{align*}
		Therefore $\|M\|_2\le 2r{L}$ and $\lambda_0^2+\lambda_1^2\le 4r^2{L}^2$.
		\item If $\lambda_1=0$ and $-\min\{{L},\lambda_{\min}(C+BA^{-1}B^\top)\}<\lambda_0\le 0$, we have
		\begin{align*}
		H(\lambda) &= C+\lambda_0 I+B \left(A+\frac{\lambda_0}{r}I\right)^{-1}B^\top \\
		& >\lambda_0 I+C+B A^{-1}B^\top\\
		&>0,
		\end{align*}
		which leads to contradiction.
	\item This part can be directly shown by combining part 2 and part 4.
	\end{enumerate}
\end{proof}
With the properties of $M$, we proceed to prove Lemma~\ref{lem:spectral_GDA_EG}.

\begin{proof}[Proof of Lemma~\ref{lem:spectral_GDA_EG}]
	Let $\lambda=\lambda_0+i\lambda_1$ be an eigenvalue of $M$, we only need to show $\rho_{1,\lambda}=|1+\etax\lambda|\le 1-\frac{1}{16r\kx}$ and $\rho_{2,\lambda}=|1+\etax\lambda+\etax^2\lambda^2|\le 1-\frac{1}{16r\kx}$, for any eigenvalue $\lambda$.

We first bound $\rho_{1,\lambda}$. There are two possible cases.
	\begin{itemize}
		\item If $\lambda_1=0$, by Lemma~\ref{lem:spectral_M}, we know $-2r{L}\le \lambda_0\le -\mu_x$. Therefore $\rho_{1,\lambda}=|1+\etax\lambda_0|\le 1-\frac{1}{4r\kappa_x}\le 1-\frac{1}{16r\kx}$.
		\item If $\lambda_1\neq 0$, we know $-2r{L}\le \lambda_0\le-\frac{\mu(r-\ky)}{2}$ and $|\lambda_1|\le \sqrt{r}{L}$. Then we have\begin{align*}
		\rho_{1,\lambda}^2=& (1+\etax\lambda_0)^2+\etax^2\lambda_1^2\\
		\le &\left(1+\frac{-\mu(r-\ky)}{8rL}\right)^2+\frac{1}{16r}\\
		\le &(1-\frac{1}{16\ky})^2+\frac{1}{32\ky}\le 1-\frac{1}{16\ky}.
		\end{align*}
		Therefore \begin{align*}
		\rho_{1,\lambda}<\sqrt{1-\frac{1}{16\kappa}}\le 1-\frac{1}{32\kappa}\le 1-\frac{1}{16r\kappa_x} .
		\end{align*}
	\end{itemize}

Next, we bound $\rho_{2,\lambda}$. Similarly, there are also two possible cases.
    \begin{itemize}
        \item If $\lambda_1=0$, noting that $-2rL\le\lambda_0\le-\mux$ and $\abs{\etax\lambda_0}\le 1/2$, we have
        \begin{align*}
\rho_{2,\lambda}\triangleq \abs{1+\etax\lambda_0+\etax^2\lambda_0^2} \le 1+\frac{1}{2}\etax\lambda_0\le 1-\frac{1}{8r\kx}.
        \end{align*}
\item If $\lambda_1\neq 0$,
\begin{align*}
    \rho_{2,\lambda}^2 =& \abs{1+\etax\lambda+\etax^2\lambda^2}^2\\
    =& \abs{1+\etax\lambda_0+\eta_x^2(\lambda_0^2-\lambda_1^2)+i\etax\lambda_1(1+2\etax\lambda_0)}^2\\
    =&\abs{1+\etax\lambda_0+\eta_x^2(\lambda_0^2-\lambda_1^2)}^2+\etax^2\lambda_1^2(1+2\etax\lambda_0)^2\\
    \le & \left(1+\frac{1}{2}\etax\lambda_0\right)^2+\etax^2\lambda_1^2\\
    \le & \left(1-\frac{\mu}{4}\etax(r-\ky)\right)^2+\frac{1}{16r}\\
    \le & 1-\frac{1}{64\kappa}.
\end{align*}
Therefore \begin{align*}
		\rho_{2,\lambda}<\sqrt{1-\frac{1}{64\kappa}}\le 1-\frac{1}{128\kappa}\le 1-\frac{1}{64r\kappa_x} .
		\end{align*}
    \end{itemize}
Noting that $\rho_{1}=\max_{\lambda}\rho_{1,\lambda}$ and $\rho_{2}=\max_{\lambda}\rho_{2,\lambda}$, we complete the proof.
\end{proof}

\section{Proofs of the convergence results for quadratic NC-SC functions}
\label{app:quad_conv}

\subsection{Proof of Theorem~\ref{thm:converence_nond_1}}

In this section, we show the linear convergence results for quadratic functions. Given Lemma~\ref{lem:spectral_GDA_EG}, if $M$ is diagonalizable, it is straight-forward to obtain a linear convergence rate. However, in some rare cases, $M$ may not be diagonalizable and things become trickier. Therefore, we further need the following lemma.
\begin{lemma}
	Let $A\in \C^{n\times n}$ be a (possibly nondiagonalizable) matrix and its spectral radius is $\rho<1$. Suppose its Jordan decomposition is $A=P^{-1}JP$, then we have\begin{align*}
	\|A^k\|_2\le sC_P k^{s-1} \rho^{k-s+1},
	\end{align*}
	where $s$ is the size of the largest Jordan block in $J$ and $C_P=\|P^{-1}\|_2\|P\|_2$ is the condition number of $P$.
	\label{lem:power_conv}
\end{lemma}

\begin{proof}[Proof of Lemma~\ref{lem:power_conv}]
	Note that
	\begin{align*}
	\norm{A^k}_2 = \norm{P^{-1}J^kP}_2\le C_P \norm{J^k}_2.
	\end{align*}
	
	So it suffices to bound $\norm{J^k}_2$. Use $J(\lambda,m)$ to denote the Jordan block with eigenvalue $\lambda$ and size $m$. We know $[J(\lambda,m)]^k=(a_{i,j})$ where $a_{i,j}={k \choose j-i}\lambda^{k-(j-i)}$ (here we define ${N\choose n}=0$ if $n<0$ or $n>N$). Then since $|\lambda|\le \rho <1$, we can bound $\|[J(\lambda,m)]^k\|_2\le m k^{m-1} \rho^{k-m+1} \le s k^{s-1} \rho^{k-s+1} $. Then we have $\|J^k\|_2\le s k^{s-1} \rho^{k-s+1}$ and complete the proof.
\end{proof}
With Lemma~\ref{lem:power_conv}, it is straight forward to show Theorem~\ref{thm:converence_nond_1} and we omit the proof.

\subsection{SGDA on quadratic functions}
Now we show Theorem~\ref{thm:conv_general_sgda}, the converge result for SGDA on quadratic functions.

\begin{proof}[Proof of Theorem~\ref{thm:conv_general_sgda}]
	We first prove the convergence for GDA.	Denote $\tz=\vz-\vz^\ast$ and $\vxi^k = (\vxi_{\vx}^k,\vxi_{\vy}^k)$ where $\vxi_{\vx}^k,\vxi_{\vy}^k$ are defined in \eqref{eq:gradient_noise}. The dynamics of GDA is
	\begin{align*}
	\tz^{k+1}=\left(I+\etax M\right)\tz^k+\etax\vxi^k,
	\end{align*}
where we can bound $\E\left[\norm{\vxi^k}_2^2\right]\le 2r^2\sigma^2/S$.
 Note that we have 
	\begin{align*}
	\left\|\prod_{k\in S} \left(I+\etax M\right)\right\|_2\le& C_P \left(1-\frac{1}{64r\kx}\right)^{|S|}.
	\end{align*} 
	Then we can bound
	\begin{align*}
	\E\left[\norm{\tz^T}_2^2\right]\le& \left\| \left(I+\etax M\right)^T\right\|_2^2 \norm{\tz^0}_2^2+\etax^2\sum_{i=0}^{T-1}\left\|\left(I+\etax M\right)^{T-i-1}\right\|_2^2 \E[\norm{\vxi^i}_2^2]\\
	\le&C_P^2\left(1-\frac{1}{64r\kx}\right)^{2T}\norm{\tz^0}_2^2+\etax^2 C_P^2 \frac{2r^2\sigma^2}{S}  \sum_{i=0}^{T-1}\left(1-\frac{1}{64r\kx}\right)^{T-i-1}\\
	\le& C_P^2\left(1-\frac{1}{64r\kx}\right)^{2T}\norm{\tz^0}_2^2+\frac{8r\kx C_P^2\sigma^2}{L^2 S},
	\end{align*}
	where the first inequality is due to the independence of $\{\vxi^k\}_{k\ge 0}$. If we require $\E\left[\norm{\tz^T}_2^2\right]\le\epsilon^2$, we only need to have
	\begin{align*}
	    T=&O\left(r\kx\log\left(\frac{C_P\norm{\tz^0}_2}{\epsilon^2}\right)\right),\\
	    S=&O\left(\frac{r\kx C_P^2\sigma^2}{L^2\epsilon^2}\right).
	\end{align*}
	This implies a gradient complexity bound of 
	\begin{align*}
	    T\cdot S=O\left(\frac{r^2\kx^2 C_P^2\sigma^2}{L^2\epsilon^2}\log\left(\frac{C_P\norm{\tz^0}_2}{\epsilon^2}\right)\right).
	\end{align*}
	Now we prove the convergence result for EG. Note that the dynamics of EG is
	\begin{align*}
	    \tz^{k+1}=&\tz^{k}+\etax M \tz^{k+1/2}+\etax\vxi^{k+1/2}\\
	    =&\tz^{k}+\etax M \left(\tz^{k}+\etax M \tz^{k}+\etax\vxi^{k}\right)+\etax\vxi^{k+1/2}\\
	    =&(I+\etax M+\etax^2M^2)\tz^{k}+\etax\vxi^{k+1/2}+\etax^2M\vxi^{k}.
	\end{align*}
	Since $\etax\norm{M}_2\le 1$ by Lemma~\ref{lem:spectral_M}, we can bound $\E\left[\norm{\vxi^{k+1/2}+\etax M\vxi^{k}}_2^2\right]\le 4r^2\sigma^2/S$. Also note that by Lemma~\ref{lem:spectral_GDA_EG}, the transition matrix of EG also satisfies $\norm{I+\etax M+\etax^2 M^2}_2\le 1-\frac{1}{64r\kx}$.
	Then the rest of proof is the same as that for GDA.
\end{proof}

\section{Proof of Theorem~\ref{thm:main_local_conv}}
\label{app:main_local_conv}
In this section, we prove our main theorem (Theorem~\ref{thm:main_local_conv}). Since the function $f$ is locally approximately quadratic, we are essentially extending the results of quadratic functions to approximately quadratic functions. So the proofs in this section build upon the lemmas and proofs for quadratic functions in Appendix~\ref{app:spectral} and Appendix~\ref{app:quad_conv}.

Note that since the Hessian of $f$ is continuous under our assumption, for any $\Delta>0$, there exist some $\delta_m>0$ such that if $\norm{\vz-\vz^\ast}_2\le \delta_m$, we have
\begin{align*}
    \max\left\{  \|\nabla^2_{\vx\vx}f(\vx;\vy)-C\|_2,\|\nabla^2_{\vx\vy}f(\vx;\vy)-B\|_2,
    \|\nabla^2_{\vy\vy}f(\vx;\vy)+A\|_2
    \right\}\le \Delta.
\end{align*}

 Before proving the theorem, we first present the following useful lemma.
\begin{lemma}
	For any two points $\vz_1$ and $\vz_2$ such that $\max\{\norm{\vz_1-\vz^\ast}_2,\norm{\vz_2-\vz^\ast}_2\}\le \delta_m$, there exists a symmetric matrix $\Tilde{H}=\begin{pmatrix}
	\Tilde{C}&\Tilde{B}\\
	\Tilde{B}^\top &-\Tilde{A}
	\end{pmatrix}$ such that
	\begin{align*}
		\nabla f(\vz_2)-\nabla f(\vz_1)= \Tilde{H}\cdot (\vz_2-\vz_1),
	\end{align*}
	where we also have
\begin{align*}
    \max\left\{\norm{C-\Tilde{C}}_2, \norm{B-\Tilde{B}}_2, \norm{A-\Tilde{A}}_2\right\}\le \Delta.
\end{align*}
	\label{lem:mvt}
\end{lemma}
\begin{proof}[Proof of Lemma~\ref{lem:mvt}]
	For any $0\le t\le 1$, denote $\vz_t=(1-t)\vz_1+t\vz_2$. We can denote the hessian at $\vz_t$ as
	\begin{align*}
		\nabla^2 f(\vz_t)=\begin{pmatrix}
		C_t&B_t\\B_t^\top & -A_t
		\end{pmatrix}.
	\end{align*}
	By Jensen's inequality, we know $\norm{\vz_t-\vz^\ast}_2\le \delta_m$. Therefore
\begin{align}
    \max\{\norm{C-C_t}_2, \norm{B-B_t}_2, \norm{A-A_t}_2\}\le \Delta.
    \label{eq:abc}
\end{align}
	Note that
	\begin{align*}
		\nabla f(\vz_2)-\nabla f (\vz_1) = \int_{0}^{1}\nabla^2 f(\vz_t)\cdot  (\vz_2-\vz_1)\,dt=:\Tilde{H}\cdot (\vz_2-\vz_1),
	\end{align*}
	where we define
	\begin{align*}
	\Tilde{H}:=\int_{0}^{1}\nabla^2 f(\vz_t)\,dt=\E_{t\sim\cU(0,1)}\left[\nabla^2 f(\vz_t)\right],
	\end{align*}
	where $\cU(0,1)$ denotes the uniform distribution over $[0,1]$.
Then we have
	\begin{align*}
		\Tilde{A}=\E_{t\sim\cU(0,1)}[A_t],\quad \Tilde{B}=\E_{t\sim\cU(0,1)}[B_t],\quad \Tilde{C}=\E_{t\sim\cU(0,1)}[C_t].
	\end{align*}
	Taking expectation of the \eqref{eq:abc} and applying Jensen's inequality, we complete the proof.
\end{proof}

With Lemma~\ref{lem:mvt}, we are ready to prove Theorem~\ref{thm:main_local_conv}.	
\begin{proof}[Proof of Theorem~\ref{thm:main_local_conv}]
Since $\nabla f(\vz^\ast)=\0$, by Lemma~\ref{lem:mvt}, for any point $\vz$ such that $\norm{\vz-\vz^\ast}_2\le\delta_m$, there exists a matrix  $\Tilde{H}=\begin{pmatrix}
	\Tilde{C}&\Tilde{B}\\
	\Tilde{B}^\top &-\Tilde{A}
	\end{pmatrix}$ such that
	\begin{align*}
	\nabla f(\vz) = 
	\Tilde{H}\cdot (\vz-\vz^\ast).
	\end{align*}
We will use $\Tilde{H}_k=\begin{pmatrix}
	\Tilde{C}_k&\Tilde{B}_k\\
	\Tilde{B}_k^\top &-\Tilde{A}_k
	\end{pmatrix}$ to such a matrix for $\vz^k$, i.e., $\nabla f(\vz^k) = 
	\Tilde{H}_k\cdot (\vz^k-\vz^\ast)$.
	
	Let use first prove the convergence rate of GDA. Choosing $\Delta\le \frac{c_0L}{r\kx C_P}$ for some small enough numerical constant $c_0$ and initialize $\vz^0$ such that $\norm{\vz^0-\vz^\ast}_2\le \delta_m/C_P$. We claim that $\norm{\vz^t-\vz^\ast}_2\le\delta_m$ for all $t\ge 0$. We will prove this claim by induction. First note that it is trivially true when $t=0$. Let us assume it is true for any $t\le T$ and try to show it also holds for $t=T+1$.
	
	Denote $\tz=\vz-\vz^\ast$. We have for any $k\le T$,
	\begin{align*}
	\tz^{k+1} =\left(I+\etax \Tilde{H}_k \right)\tz^{k}=\left(I+\etax M+E_k\right)\tz^{k},
	\end{align*}
	where 
	\begin{align*}
	E_k:=\etax(\Tilde{H}_k-M)=\etax\begin{pmatrix}
	-\left(\Tilde{C}_k-C \right) & -\left(\Tilde{B}_k-B \right) \\ r \left(\Tilde{B}_k-B \right)^\top & r \left(\Tilde{A}_k-A \right)
	\end{pmatrix}.
	\end{align*}
	It is easy to bound $\|E_k\|_2 \le 2\etax  r\Delta$.
	Suppose $I+\etax M=P^{-1}\Lambda P$ is diagonalizable, where $\Lambda=\text{diag}(\lambda_1,\ldots,\lambda_{m+n})$ and we have shown $|\lambda_i|\le 1-\frac{1}{64r\kx}$. Then we can bound that
	\begin{align*}
	\left\|\prod_{k=0}^T \left(I+\etax M+E_k\right)\right\|_2\le& \sum_{i=0}^{T+1} {T+1\choose i} C_p^{i+1}\left(1-\frac{1}{64r\kx}\right)^{T+1-i}\left(2\etax  r\Delta\right)^i\\
	\le & C_P \left(1-\frac{1}{64r\kx}+2\etax  r\Delta C_P\right)^{T+1}\\
	\le &C_P \left(1-\frac{1}{128r\kx}\right)^{T+1}.
	\end{align*} 
Then we have
\begin{align*}
    \norm{\vz^{T+1}-\vz^\ast}_2\le C_P \left(1-\frac{1}{128r\kx}\right)^{T+1}\norm{\vz^0-\vz^\ast}_2\le C_P\cdot\frac{\delta_m}{C_P}= \delta_m.
\end{align*}
Therefore we have proved the claim as well as obtained the convergence rate for GDA. 

Now we prove the convergence rate for EG. Note that the dynamics of EG can be written as
\begin{align*}
    \tz^{k+1}=&\tz^k+\etax\Tilde{H}_{k+1/2}\cdot \tz^{k+1/2}\\
    =&\left(1+\etax\Tilde{H}_{k+1/2}+\etax^2\Tilde{H}_{k+1/2}\Tilde{H}_{k}\right)\tz^k\\
    =:&\left(1+\etax M+\etax^2 M^2+F_k\right)\tz^k,
\end{align*}
where we have
\begin{align*}
    \norm{F_k}_2 =\norm{(1+\etax M)E_{k+1/2}+E_k(\etax M+E_{k+1/2})}_2 \le 4\norm{E_k}_2,
\end{align*}
where in the inequality we use the fact $\norm{E_{k+1/2}}_2\le 1$ and $\etax\norm{M}_2\le 1$ due to Lemma~\ref{lem:spectral_M}. Also note that by Lemma~\ref{lem:spectral_GDA_EG}, the transition matrix of EG also satisfies $\norm{I+\etax M+\etax^2 M^2}_2\le 1-\frac{1}{64r\kx}$. The rest of proof for EG is similar to that of GDA.

\end{proof}
\section{Proof of Theorem~\ref{thm:main_local_conv_stochastic}}
\label{app:main_local_conv_stochastic}

\begin{proof}[Proof of Theorem~\ref{thm:main_local_conv_stochastic}]
	We use the same notation as in the proof of Theorem~\ref{thm:conv_general_sgda}. We first prove the convergence for GDA.
	Similar to the proof of Theorem~\ref{thm:main_local_conv}, we choose $\Delta\le \frac{c_0L}{r\kx C_P}$ for some small enough numerical constant $c_0$ and initialize $\vz^0$ such that $\norm{\vz^0-\vz^\ast}_2\le \delta_m/2C_P$. We claim that with probability at least $1-p$ where $p$ is an arbitrarily small constant, we have $\norm{\vz^t-\vz^\ast}_2\le\delta_m$ for any $t\le T$ where $T$ is the maximum number of iterations. Denote the event
	\begin{align*}
	    \cE = \left\{\vxi\Big|\etax C_P\sum_{k=0}^{T-1}\|\vxi^k\|_2\le \epsilon/2\right\},
	\end{align*}
	where $\vxi=\{\vxi^k\}_{k=0}^{T-1}$ and we can bound $\E\left[\norm{\vxi^k}_2^2\right]\le 2r^2\sigma^2/S$. Choosing
	\begin{align*}
	    S=&\frac{C T^2 C_P^2\sigma^2}{L^2\epsilon^2}
	\end{align*}
	for some large enough constant $C$,	by Markov inequality, we have
	\begin{align*}
	    \mathbb{P}\left(\vxi\notin\cE\right)\le& \frac{2\etax C_P}{\epsilon}\sum_{k=0}^{T-1}\E[\|\vxi^k\|_2]\\
	    \le & \frac{2 C_P}{4rL\epsilon}\cdot T\cdot \sqrt{2}r\sigma /\sqrt{S}\\
	    \le & p.
	\end{align*}
 In other words, with probability at least $1-p$, we have $\vxi\in\cE$. We will condition on this event in the following proof.
	
Then we can use the induction technique similar to the proof of Theorem~\ref{thm:main_local_conv}.	Suppose that $\|\tz^k\|\le\delta_m$ for all $k\le K$ where $K<T$. The dynamics of GDA is
	\begin{align*}
	\tz^{k+1}=\left(I+\etax M+E_k\right)\tz^k+\etax\vxi^k,
	\end{align*}
where $E_k$ is defined in the proof of Theorem~\ref{thm:main_local_conv}.
 Note that in this case, following the proof of Theorem~\ref{thm:main_local_conv}, we can bound
	\begin{align*}
	\left\|\prod_{k\in \cI} \left(I+\etax M+E_k\right)\right\|_2\le& C_P \left(1-\frac{1}{128r\kx}\right)^{|\cI|},\text{ for any index set }\cI.
	\end{align*} 
The problem is more difficult than that in Theorem~\ref{thm:conv_general_sgda} because $E_k$ could potentially depend on the gradient noise $\{\vxi^i\}_{i<k}$. Conditioned on the event $\vxi\in\cE$, we can bound that
	\begin{align*}
	\norm{\tz^{K+1}}_2\le& \left\| \prod_{k=0}^K\left(I+\etax M+E_k\right)\right\|_2 \norm{\tz^0}_2+\etax\sum_{i=0}^{K}\left\|\prod_{k=i}^K\left(I+\etax M+E_k\right)\right\|_2 \norm{\vxi^i}_2\\
	\le&C_P\left(1-\frac{1}{128r\kx}\right)^{K+1}\norm{\tz^0}_2+\etax C_P  \sum_{i=0}^{T-1} \norm{\vxi^i}_2\\
	\le& C_P\norm{\tz^0}_2+\epsilon/2\le \delta_m,
	\end{align*}
		where in the last inequality we use $\epsilon\le\delta_m$ which holds because otherwise the initialization point is already $\epsilon$-approximately optimal.
	Therefore the induction hypothesis also holds for $K+1$. Then we can conclude that $\|\tz^k\|_2\le\delta_m$ for all $0\le k\le T-1$, i.e., the iterates never escape from the neighborhood around $\vz^*$ under the event $\vxi\in\cE$.
	
	From the above inequalities, when $K=T-1$, we get
	\begin{align*}
	    \norm{\tz^{T}}_2
	\le&C_P\left(1-\frac{1}{32r\kx}\right)^{T}\norm{\tz^0}_2+\epsilon/2.
	\end{align*}
	 If we require $\norm{\tz^T}_2\le\epsilon$, we only need to have
	\begin{align*}
	    T=&O\left(r\kx\log\left(\frac{C_P\norm{\tz^0}_2}{\epsilon}\right)\right),
	\end{align*}
	which implies
	\begin{align*}
	    S=&O\left(\frac{r^2\kx^2 C_P^2\sigma^2}{L^2\epsilon^2}\log^2\left(\frac{C_P\norm{\tz^0}_2}{\epsilon}\right)\right).
	\end{align*}
	This implies a gradient complexity bound of 
	\begin{align*}
	    T\cdot S=O\left(\frac{r^3\kx^3 C_P^2\sigma^2}{L^2\epsilon^2}\log^3\left(\frac{C_P\norm{\tz^0}_2}{\epsilon^2}\right)\right).
	\end{align*}
	Now we prove the convergence result for EG. Note that the dynamics of EG is
	\begin{align*}
	    \tz^{k+1}=&\tz^{k}+(\etax M+E_{k+1/2}) \tz^{k+1/2}+\etax\vxi^{k+1/2}\\
	    =&\tz^{k}+(\etax M+E_{k+1/2}) \left(\tz^{k}+(\etax M+E_{k}) \tz^{k}+\etax\vxi^{k}\right)+\etax\vxi^{k+1/2}\\
	    =&(I+\etax M+\etax^2M^2+F_k)\tz^{k}+\etax\vxi^{k+1/2}+\etax(\etax M+E_{k+1/2})\vxi^{k}.
	\end{align*}
	Note that it is straight forward to bound $\norm{\etax M+E_{k+1/2}}_2\le 1$.
Then we can bound $\E\left[\norm{\vxi^{k+1/2}+(\etax M+E_{k+1/2})\vxi^{k}}_2\right]\le 4r\sigma/\sqrt{S}$. Then the rest of proof is similar to that for GDA. In particular, we can define the event 
	\begin{align*}
	    \cE = \left\{\vxi\Big|\etax C_P\sum_{k=0}^{T-1}(\|\vxi^k\|_2+\|\vxi^{k+1/2}\|_2)\le \epsilon/2\right\}
	\end{align*}

and show it holds with high probability if choosing a large enough $S$. Following the proof of Theorem~\ref{thm:main_local_conv} for EG, we can also show that 
	\begin{align*}
	\left\|\prod_{k\in \cI} \left(I+\etax M+\etax^2M^2+F_k\right)\right\|_2\le& C_P \left(1-\frac{c_0}{r\kx}\right)^{|\cI|},\text{ for any index set }\cI,
	\end{align*} 
	for some small enough numerical constant $c_0$. Then the the following steps are straight forward.

\end{proof}
\section{Additional experimental results}

\begin{figure*}[t]
     \centering
     \begin{subfigure}[b]{0.48\textwidth}
         \centering
         \includegraphics[width=\textwidth]{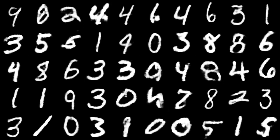}
         
         \vspace{0.2cm}
         \includegraphics[width=\textwidth]{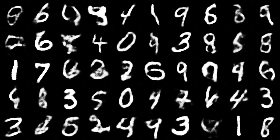}
         
         \vspace{0.2cm}
         \includegraphics[width=\textwidth]{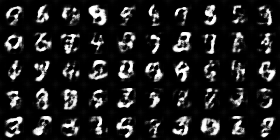}
         \caption{}
         \label{fig:images_app_a}
     \end{subfigure}
     \hfill
     \begin{subfigure}[b]{0.48\textwidth}
         \centering
         \includegraphics[width=\textwidth]{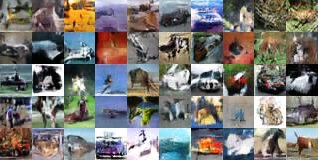}
         
         \vspace{0.2cm}
         \includegraphics[width=\textwidth]{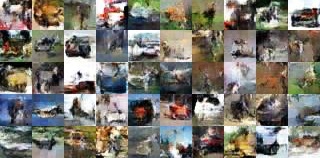}
         
         \vspace{0.2cm}
         \includegraphics[width=\textwidth]{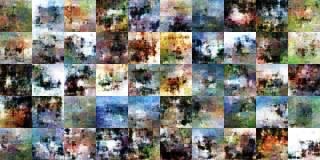}
         \caption{}
         \label{fig:images_app_b}
     \end{subfigure}
        \caption{Generated images of the convergent generators for MNIST (a) and CIFAR10 (b). All the generators are trained using the WGAN-GP model with SGDA. The stepsize choices $(\etax,\etay)$ are set to be $(0.001,0.001)$, $(0.0001,0.001)$, and $(1e-5,0.001)$ from the first row to the third row.}
        \label{fig:images_app}
\end{figure*}

\begin{figure*}[t]
     \centering
         \includegraphics[width=0.48\textwidth]{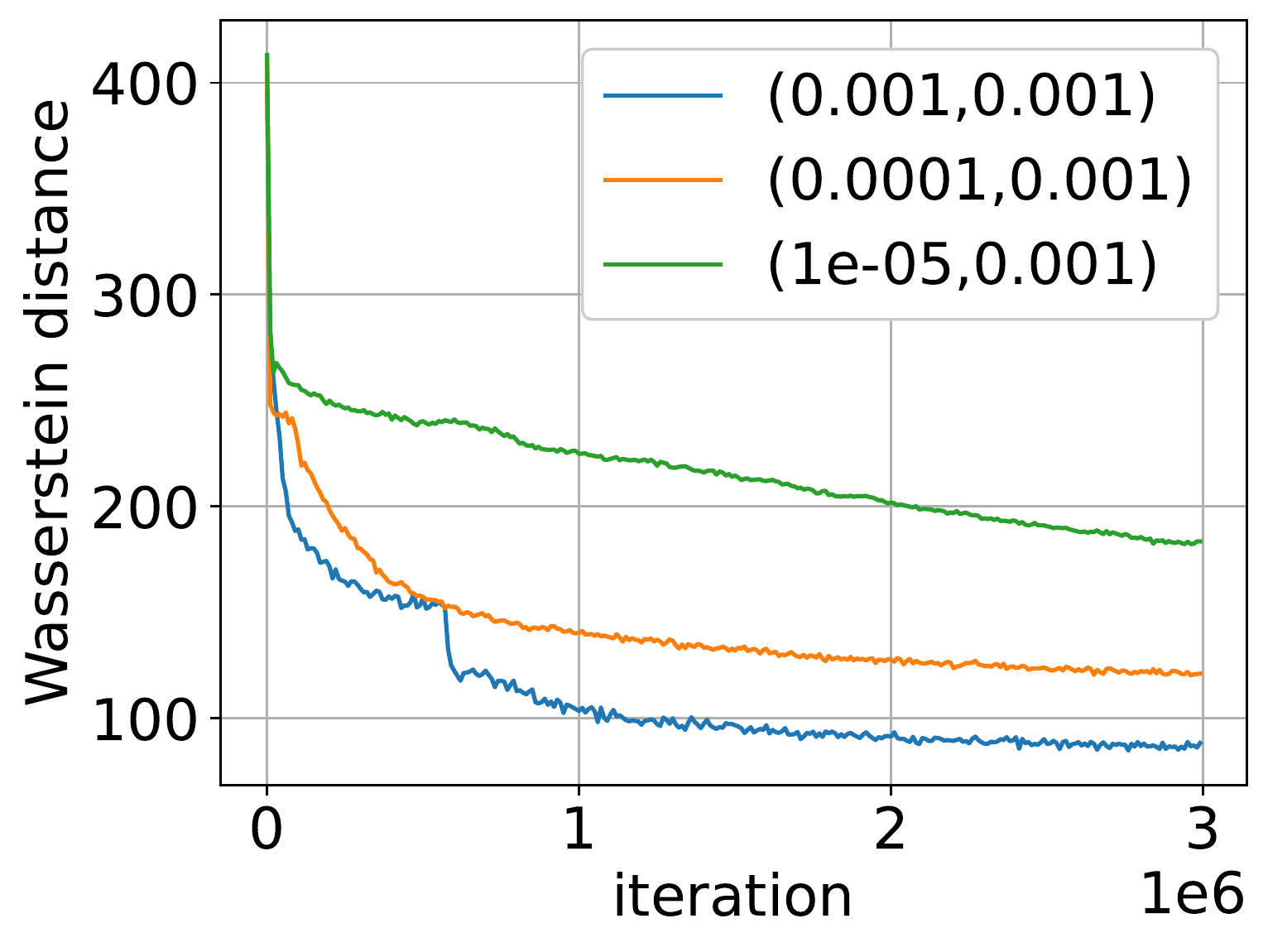}
        \caption{Evolution of FID score on CIFAR10. We train WGAN-GP models under different stepsize ratios $(\etax,\etay)=(0.001,0.001),(0.0001,0.001),(1e-05,0.001)$.}
        \label{fig:fid}
\end{figure*}

\begin{figure*}[t]
     \centering
     \begin{subfigure}[b]{0.48\textwidth}
         \centering
         \includegraphics[width=\textwidth]{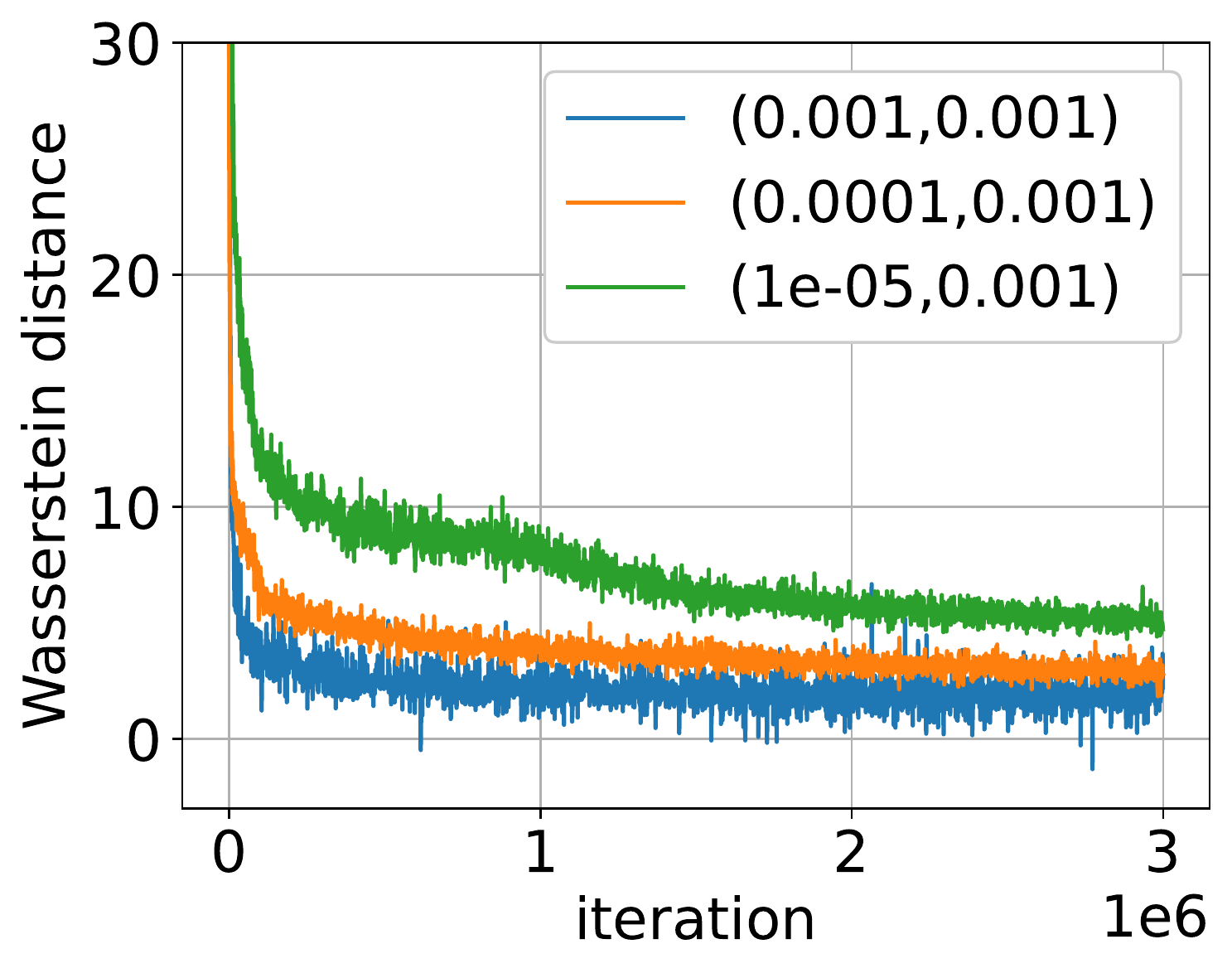}
         \caption{}
         \label{fig:dcgan_curve_a}
     \end{subfigure}
     \hfill
     \begin{subfigure}[b]{0.48\textwidth}
         \centering
         \includegraphics[width=\textwidth]{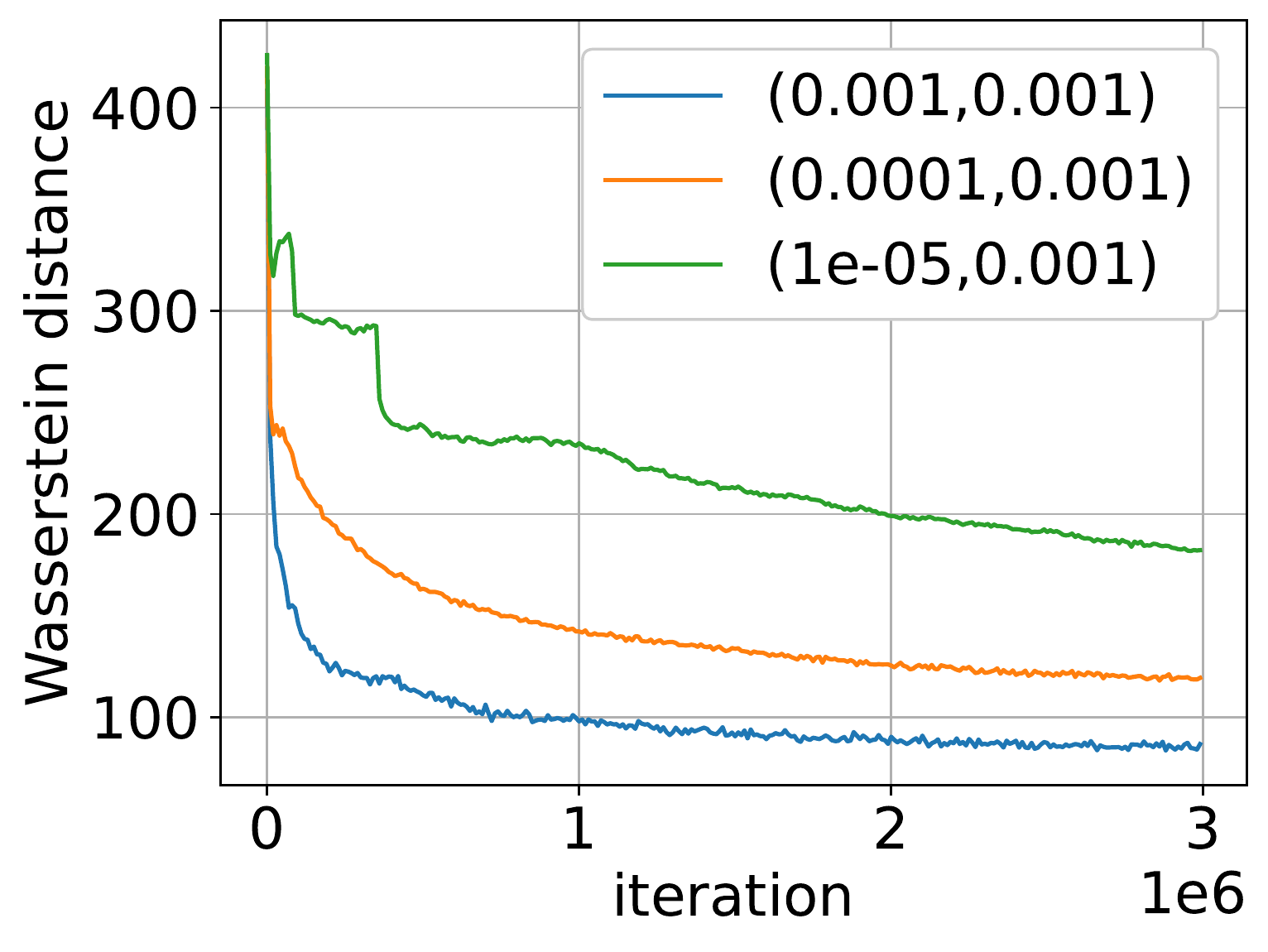}
         \caption{}
         \label{fig:dcgan_curve_b}
     \end{subfigure}
        \caption{Evolution of Wasserstein Distance (a) and FID score (b) for CIFAR10. All the generators are trained using the DCGAN with SGDA under different stepsize choices $(\etax,\etay)=$  $(0.001,0.001),(0.0001,0.001),(1e-5,0.001)$.}
        \label{fig:dcgan_curve}
\end{figure*}

\begin{figure*}[t]
     \centering

         \includegraphics[width=0.7\textwidth]{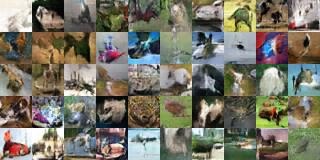}
         
         \vspace{0.2cm}
         
         \includegraphics[width=0.7\textwidth]{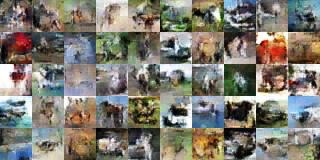}
         
         \vspace{0.2cm}
         \includegraphics[width=0.7\textwidth]{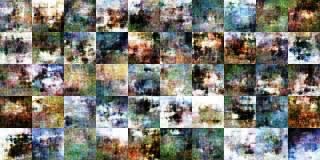}

        \caption{Generated images of the convergent generators for DCGAN and CIFAR10. The stepsize choices $(\etax,\etay)$ are set to be $(0.001,0.001)$, $(0.0001,0.001)$, and $(1e-5,0.001)$ from the first row to the third row.}
        \label{fig:dcgan_images}
\end{figure*}

We present additional experimental results in this section. 
In Figure~\ref{fig:mot_exp}, we showed the training curves of SGDA on MNIST and CIFAR10 with different stepsize ratios. 
Now we show the generated images of the trained generator in Figure~\ref{fig:images_app}. As we can clearly see, a larger stepsize ratio not only makes it converging more slowly, but also results in a worse generator for both MNIST and CIFAR10 datasets. This is also shown in by the evolution FID score during training for CIFAR10 in Figure~\ref{fig:fid}. 
All the experimental results so far use the WGAN-GP model. We will show that another GAN model, DCGAN, has similar training behaviors. Figure~\ref{fig:dcgan_curve} shows the evolution of the Wasserstein Distance and FID score under different stepsize choices for DCGAN and CIFAR10. Figure~\ref{fig:dcgan_images} shows the generated images in this setting. As we can see, the results are similar to those for WGAN-GP.

\end{document}